\documentclass[11pt]{article}
\usepackage{anysize}
\marginsize{2.0cm}{2.0cm}{2.0 cm}{2.0 cm}
\usepackage{setspace}
\usepackage{graphicx}
\usepackage{color}
\usepackage{epstopdf}
\usepackage[titletoc]{appendix}
\usepackage{graphicx}
\usepackage{latexsym,amsfonts,amsmath,amssymb}
\newtheorem{theorem}{Theorem}
\newtheorem{lemma}{Lemma}
\newtheorem{corollary}{Corollary}
\newtheorem{definition}{Definition}
\newtheorem{proposition}{Proposition}

\def\EE{\mathsf{E}}
\def\FF{\mathbb{F}}
\def\HH{\mathbb{H}}
\def\MM{\mathsf{M}}

\def\RR{\mathbb{R}}

\def\PP{\mathsf{P}}

\begin{document}

\title{Extreme occupation measures in  Markov decision processes with a cemetery}

\date{}

\author{Alexey Piunovskiy\thanks{Corresponding author. Department of Mathematical Sciences, University of
Liverpool, Liverpool, U.K.. E-mail: piunov@liv.ac.uk.}
\and Yi
Zhang \thanks{School of Mathematics, University of Birmingham, Edgbaston,
Birmingham,
B15 2TT, U.K.. Email: y.zhang.29@bham.ac.uk} }

\maketitle

\par\noindent\textbf{Abstract}.
In this paper, we consider a Markov decision process (MDP) with a Borel state space $\textbf{X}\cup\{\Delta\}$, where $\Delta$ is an absorbing state (cemetery), and a Borel action space $\textbf{A}$. We consider the space of finite occupation measures restricted on $\textbf{X}\times \textbf{A}$, and the extreme points in it. It is possible that some strategies have infinite occupation measures. Nevertheless, we prove that every finite extreme occupation measure is generated by a deterministic stationary strategy. Then, for this MDP, we consider a constrained problem with total undiscounted criteria and $J$ constraints, where the cost functions are nonnegative. By assumption, the strategies inducing infinite occupation measures are not optimal. Then, our second main result is that, under mild conditions, the solution to this constrained MDP is given by a mixture of no more than $J+1$ occupation measures generated by deterministic stationary strategies.

\bigskip
\par\noindent\textbf{Keywords.}  Markov decision process, Total cost, Occupation measure, Mathematical programming, Extreme point, Mixture
\bigskip

\par\noindent{\bf AMS 2000 subject classification:}  Primary 90C40; Secondary 90C25

\section{Introduction}\label{pz23bSec01}

Perhaps the first paper, where the discounted Markov decision process (MDP) was reformulated as a linear program, is \cite{depen}. The modern so called `convex analytic approach' originates from the works by V.Borkar \cite{Borkar:1988,borkar91}. It is applied to the models with total cost (discounted or not) as well as with the long-run average cost: let us only mention the book treatments \cite{altman,Hernandez-Lerma:1999,Piunovskiy:1997} and the survey \cite{Borkar:2002}. This approach proved to be especially fruitful in dealing with problems with constraints, see the survey \cite{Piunovskiy:1998a} and the authoritative monograph \cite{kall} on finite MDPs, i.e., MDPs with finite state and action spaces.  For convex analytic approach for continuous-time MDPs, see e.g., \cite{guozhang,guohuangzhang,APZY}.

The convex analytic approach is based on the reformulation of the constrained MDP problem as a convex optimization problem in the space of occupation measures with affine objective functions and inequality constraints, where the occupation measures are defined in accordance with the performance criteria of the MDP problem. The space of occupation measures is a convex space (i.e., a convex subset of a cone, not necessarily of a vector space). Thus, here, the relevant notions, such as convex optimization problem, affine functions and extreme points, are understood with respect to (wrt) the underlying convex space, see \cite{archive}. An important target to show is the existence of an optimal strategy for the MDP problem, whose occupation measure is the convex combination of finitely many extreme points in the space of occupation measures, which we call extreme occupation measures. If the number of constraints in the MDP problem is $J$, the mixture is over at most $J+1$ extreme occupation measures. Such a strategy is called a $(J+1)$-mixed optimal strategy. Then a key ingredient in the convex analytic approach to MDPs is the characterization of such extreme occupation measures. This task is easier when the state is discrete (finite or countable), as considered in \cite{altman,Borkar:1988,borkar91,kall}, but our consideration in this paper is a Borel MDP model, by which we mean an MDP with Borel state and action spaces.

Let us concentrate on the literature for Borel MDP models. For discounted MDPs, the most relevant recent works include \cite{Dufour:2016,FeiJN,Gonzalez:2000}, where by using the convex analytic approach, optimal stationary strategies were proved. While mixed strategies were not considered in \cite{Dufour:2016,Gonzalez:2000}, in establishing the existence of a so called optimal chattering strategy in \cite{FeiJN} (see also \cite{gonz,gonzvill}), the existence of an optimal $(J+1)$-mixed strategy was observed, see the proof of \cite[Theorem 2]{FeiJN}. For discounted MDPs but under more restrictive conditions, this result appeared in \cite{Piunovskiy:1997,Zhang:2013}. In an absorbing MDP, there is a costless cemetery in the state space, and given the initial state, under each strategy, the expected time until the state process reaches the cemetery is finite. In fact this is equivalent to the expected absorbing time being bounded in the set of all strategies, see \cite[p.132]{feiroth}. It is well known that discounted MDPs are special cases of absorbing MDPs, by viewing the discount factor as the parameter of a geometrically distributed external killing time. For absorbing MDPs, the existence of a $(J+1)$-mixed optimal strategy was obtained by Feinberg and Rothblum, see \cite[Theorem 9.2]{feiroth}, as well as that each extreme occupation measure is generated by a deterministic stationary strategy, see \cite[Lemma 4.6]{feiroth}. The convex analytic approach was also developed for optimal stopping problems \cite{dufour10}.

In the present paper, we consider an MDP with a Borel state space $\textbf{X}\cup\{\Delta\}$ and a Borel action space $\textbf{A}$. The point $\Delta$ is a single cemetery. We call such a model an MDP with a cemetery, though it is also known under other names such as the stochastic shortest path problem, see \cite{Bertsekas:2018}, where unconstrained MDP problems were considered and the main interest was the characterization of the optimal value function out of the class of so called proper strategies in terms of the solution to the optimality equation. It is without loss of generality that we consider a fixed initial state rather than a fixed initial distribution. We also assume that the cemetery $\Delta$ is costless. For this reason, we consider occupation measures as the total expected state-action frequencies restricted on $\textbf{X}\times \textbf{A}$. If a strategy has a finite occupation measure, we call it an absorbing strategy with the given initial state. Proper strategies  as considered in \cite{Bertsekas:2018} can be viewed as special types of absorbing strategies. If the occupation measure of each strategy is finite, our model becomes the absorbing model. Nevertheless, similarly to \cite{Dufour:2012}, here we do allow that some strategies have infinite occupation measures. This is the main novelty compared with the aforementioned works \cite{Dufour:2016,feiroth,FeiJN,Gonzalez:2000}, see more comments on this below.

Our contributions are as follows. First, we show that every extreme point of the space of finite occupation measures is generated by a deterministic stationary strategy. Then, we consider a constrained problem with total undiscounted criteria and $J$ constraints, where the cost functions are nonnegative.  We formulate the problem as a convex program in the space of occupation measures (see (\ref{eB7})). Under mild conditions, we show that there exists an optimal strategy whose occupation measure is in the form of a mixture of no more than $J+1$ occupation measures of deterministic stationary strategies.

For the latter result, we make the assumption, which in particular, implies that strategies inducing infinite occupation measures are not optimal. Under this assumption, for the MDP problem, instead of dealing with the whole space of occupation measures, it is sufficient to work with the space of finite occupation measures. However, restricting an MDP to absorbing strategies is not the same as considering an absorbing MDP. Firstly, in an absorbing MDP, the total values of all occupation measures are bounded above, whereas if the MDP is not absorbing, then the total values of all finite occupation measures can be unbounded. This can be seen by considering an optimal stopping problem as in \cite{dufour10}: for the set of  strategies,  stopping at step  $n=1,2,\dots$, the values of their (finite) occupation measures are unbounded. In this connection, we mention that, for discounted MDPs, see e.g., \cite{Piunovskiy:1997}, it is convenient to endow the space of occupation measures with the weak topology generated by bounded continuous functions. The same was done in \cite{feiroth} for absorbing MDPs. To deal with infinite occupation measures, we endow that space with the final topology generated by the projection mapping from the space of strategic measures to occupation measures. These features require new proofs of the key theorems on the characterization of the extreme finite occupation measures (see Theorem \ref{lB33}) and on the sufficiency of mixtures of (occupation measures of) deterministic strategies (see Theorem \ref{tB2}).

In terms of other relevant works, we mention the following. First, constrained total undiscounted Borel MDPs with non-negative cost functions were also studied in  \cite{Dufour:2012}. Although it was not assumed a priori in \cite{Dufour:2012} that there is a costless cemetery in the state space, it was shown under some conditions  that one can always construct a costless cemetery set, after modifying the admissible action spaces on that set. By merging this set as a costless cemetery, we may view the model in \cite{Dufour:2012} in the framework of the present paper, and can apply to it our first result on the characterization of extreme finite occupation measures. Except for special cases, our second result concerning the optimal mixed strategies are not applicable to the model in \cite{Dufour:2012} because no assumption was made in \cite{Dufour:2012} that strategies with infinite occupation measures were suboptimal or infeasible. On the other hand, neither the extreme occupation measures nor the mixed strategies were considered in \cite{Dufour:2012}. The paper here can be viewed as a complement to it.
Second, we note that the results in this paper are also relevant to the studies in continuous-time MDPs, see \cite{PiunovskiyZhang:2020,PiunovskiyZhang:2022SICON}, because the problems considered therein were eventually reduced to an MDP model,  see more details in the book \cite{APZY}.

Allowing the cost functions to be negative leads to a more complicated theory. The convex analytic approach to such constrained MDPs was developed in \cite{dufour13,dufour20}, but mixtures of occupation measures were not considered there.

The rest of this paper is organized as follows.  The MDP model under study is described in Section \ref{s1}. Several necessary auxiliary statements are given in Section \ref{s2}, including the known results on the solvability of the formulated problem. Sections \ref{s4} and \ref{s5} present the main results: characterization of the extreme occupation measures and sufficiency of the finite mixtures of deterministic stationary strategies. The paper ends with a conclusion in Section \ref{s6}. The proofs of the main statements are postponed to the appendix.

\section{Description of the model}\label{s1}

The primitives of an {MDP}  are the following.
\begin{itemize}
\item The {state space} is ${\bf X}_\Delta={\bf X}\cup\{\Delta\}$, where ${\bf X}$ is a nonempty topological Borel space, endowed with the $\sigma$-algebra ${\cal B}(\textbf{X})$, and $\Delta$ is the isolated absorbing state (cemetery).
\item The {action space}  ${\bf A}$ is a nonempty topological Borel space, endowed with the  $\sigma$-algebra ${\cal B}(\textbf{A}).$
\item The {transition probability}  $p(dy|x,a)$ is a stochastic kernel from ${\bf X}_\Delta\times{\bf A}$ to ${\cal B}({\bf X}_\Delta)$; $p(\{\Delta\}|\Delta,a)\equiv 1$.
\item The $[0,+\infty]$-valued one-step {cost functions}  $r_j(\cdot,\cdot)$ on $\textbf{X}_\Delta\times \textbf{A}$, $j=0,1,\ldots,J$, where $J\in\{0,1,\dots\}$ is a fixed integer; $r_j(\Delta,a)\equiv 0$.
\end{itemize}

Usually, the initial state $x_0\in{\bf X}$ is fixed, but sometimes we consider other arbitrarily fixed initial states $x\in{\bf X}$. (See, e.g., Lemma \ref{lB9}.)

Regarding terminology, we often refer to $\{\textbf{X}_\Delta,\textbf{A},p,\{r_j\}_{j=0}^J\}$ as a MDP model or simply a MDP. We may also consider the `cost-free' MDP model $\{\textbf{X}_\Delta,\textbf{A},p\}$ because several definitions and properties presented below do not involve the properties of the cost functions.

\begin{definition}[Strategy]\label{27MarchImBookAppendixDef01} Consider the MDP $\{\textbf{X}_\Delta,\textbf{A},p\}$.
\begin{itemize}
\item[(a)]A {strategy}   $\sigma=\{\sigma_n\}_{n=1}^\infty$ is a sequence of stochastic kernels such that for each $n=1,2,\dots,$ $\sigma_n(da|x_0,a_1,\dots,x_{n-1})$ is a stochastic kernel from $(\textbf{X}_\Delta\times \textbf{A})^{n-1}\times \textbf{X}_\Delta$ to ${\cal B}(\textbf{A})$, where $(\textbf{X}_\Delta\times \textbf{A})^{0}\times \textbf{X}_\Delta:=\textbf{X}_\Delta.$
\item[(b)]A strategy $\sigma=\{\sigma_n\}_{n=1}^\infty$ is {Markov}  if for each $n=1,2,\dots,$ there is a stochastic kernel $\sigma^M_n(da|x_{n-1})$ from $\textbf{X}_\Delta$ to ${\cal B}(\textbf{A})$ such that
\begin{eqnarray*}
\sigma^M_n(da|x_{n-1})=\sigma_n(da|x_0,a_1,\dots,x_{n-1})
\end{eqnarray*}
for each $(x_0,a_1,\dots,x_{n-1})\in (\textbf{X}_\Delta\times \textbf{A})^{n-1}\times \textbf{X}_\Delta.$
\item[(c)] A strategy $\sigma=\{\sigma_n\}_{n=1}^\infty$ is called {stationary}  if there is a stochastic kernel $\sigma^s(da|x)$ from $\textbf{X}_\Delta$ to ${\cal B}(\textbf{A})$ such that
\begin{eqnarray*}
\sigma^s(da|x_{n-1})=\sigma_n(da|x_0,a_1,\dots,x_{n-1})
\end{eqnarray*}
for each $n=1,2,\dots$, and $(x_0,a_1,\dots,x_{n-1})\in (\textbf{X}_\Delta\times \textbf{A})^{n-1}\times \textbf{X}_\Delta.$ Below,  a stationary strategy is usually identified with $\sigma^s$.
\item[(d)]
If $\sigma^s(da|x)$ is concentrated on  $\{\varphi(x_{n-1})\}$,   where  $\varphi$ is an $\bf A$-valued measurable mapping, then the stationary strategy is called {deterministic stationary}. With conventional abuse of notations, we often signify a deterministic stationary strategy by $\varphi$.
\item[(e)] We always assume that $\sigma_n(\{\hat a\}|x_0,a_1,\ldots,\Delta)=1$ whenever $x_{n-1}=\Delta$. Here $\hat a\in{\bf A}$ is an arbitrarily fixed action.
\end{itemize}
\end{definition}

As is well known, for each control strategy $\sigma$ and initial state $x_0\in{\bf X}$, there is a unique {strategic measure} on the sample space $\Omega:=({\bf X}_\Delta\times{\bf A})^\infty$, denoted as ${\mathsf P}^\sigma_{x_0}$, which  is specified by the following conditions:
\begin{eqnarray}\label{ea3}
\mathsf{P}_{x_0}^\sigma(X_0\in dy)=\delta_{x_0}(dy);
\end{eqnarray}
and for each $n=1,2,\dots,$  $\Gamma_i^\textbf{X}\in {\cal B}(\textbf{X})$ ($i=0,1,\dots,n$) and $\Gamma_i^\textbf{A}\in {\cal B}(\textbf{A})$ ($i=1,2,\dots,n$),
\begin{eqnarray}\label{ea4}
&&\mathsf{P}_{x_0}^\sigma(X_0\in \Gamma_0^\textbf{X},~A_1\in \Gamma_1^\textbf{A},~\dots,~X_{n-1}\in \Gamma_{n-1}^\textbf{X},~A_n\in \Gamma_n^\textbf{A})\\
&=&\int\limits_{\Gamma_{0}^\textbf{X}\times \Gamma_1^\textbf{A}\times\dots\times \Gamma_{n-1}^\textbf{X}}\sigma_n(\Gamma_n^\textbf{A}|x_0,a_1,\dots,x_{n-1})
\mathsf{P}_{x_0}^\sigma(X_0\in dx_0,A_1\in da_1,\dots,X_{n-1}\in dx_{n-1});\nonumber
\end{eqnarray}
and
\begin{eqnarray}\label{ea8}
&&\mathsf{P}_{x_0}^\sigma(X_0\in \Gamma_0^{\textbf{X}},~A_1\in \Gamma_1^{\textbf{A}},\dots,~X_{n}\in \Gamma_{n}^{\textbf{X}})\\
&=&\int\limits_{\Gamma_{0}^\textbf{X}\times \Gamma_1^\textbf{A}\times\dots\times \Gamma_{n-1}^\textbf{X}\times \Gamma_n^\textbf{A}}p(\Gamma_n^\textbf{X}|x_{n-1},a_{n})\nonumber \\
&&\times\mathsf{P}_{x_0}^\sigma(X_0\in dx_0,A_1\in da_1,\dots,X_{n-1}\in dx_{n-1},A_{n}\in da_{n}).\nonumber
\end{eqnarray}
For details, see \cite{Dynkin:1979,Hernandez-Lerma:1996,Piunovskiy:1997}. Denote by $\Sigma$ the set  of all strategies, and by ${\cal P}:=\{\mathsf{P}^\sigma_{x_0}:~\sigma\in \Sigma\}$ the set of all strategic measures (with the initial state $x_0\in{\bf X}$).
The expectation taken with respect to $\mathsf{P}_{x_0}^\sigma$ is denoted by $\mathsf{E}_{x_0}^\sigma.$ We equip the space of probability measures on ${\cal B}(\Omega)$, denoted as ${\cal P}(\Omega)$, with the weak topology generated by bounded continuous functions on ${\Omega}$, and fix its trace $\tau$ on the space $\cal P$ of all strategic measures. Then ${\cal P}(\Omega)$ is a Borel space, see \cite[Corollary 7.25.1]{Bertsekas:1978}, and we endow ${\cal P}(\Omega)$ with its Borel $\sigma$-algebra.

The {constrained} optimal control problem for the MDP model $\{\textbf{X}_\Delta,\textbf{A},p,\{r_j\}_{j=0}^J\}$ is
\begin{eqnarray}\label{ea6}
&\mbox{Minimize over all strategies $\sigma$:~}& \mathsf{E}_{x_0}^\sigma\left[\sum_{n=0}^\infty r_0(X_n,A_{n+1})\right]\\
&\mbox{subject to}& \mathsf{E}_{x_0}^\sigma\left[\sum_{n=0}^\infty r_j(X_n,A_{n+1})\right]\le d_j,~~j=1,2,\ldots, J. \nonumber
\end{eqnarray}
If $J=0$, then the problem is called {unconstrained}.

\begin{definition}[Feasible and optimal strategies]
A strategy is called feasible if all the constraints in (\ref{ea6}) are satisfied; it is called feasible with a finite value if, additionally, $\mathsf{E}_{x_0}^\sigma\left[\sum_{n=0}^\infty r_0(X_n,A_{n+1})\right]<\infty$;
it is called optimal if it solves problem (\ref{ea6}).
\end{definition}

\begin{definition}[Semicontinuous MDP] \label{da12}
A MDP $\{\textbf{X}_\Delta,\textbf{A},p,\{r_j\}_{j=0}^J\}$  is called  {semicontinuous}  if
\begin{itemize}
\item[(a)] The action space $\textbf{A}$ is compact.
\item[(b)] For each bounded continuous function $f(\cdot)$ on $\textbf{X}$, $\int_\textbf{X} f(y)p(dy|x,a)$ is continuous in $(x,a)\in \textbf{X}\times \textbf{A}.$
\item[(c)] For each $j=0,1,\dots,J,$ the function $r_j(\cdot,\cdot)$ is lower semicontinuous on $\textbf{X}\times \textbf{A}.$
\end{itemize}
\end{definition}

\section{Preliminaries}\label{s2}
In this section, we collect some preliminary results, which will be needed in proving the main results of this paper. Several of them are known, or follow from well known results. They will be called propositions. We thus skip the proofs of most of them but always refer to relevant literature.

\begin{proposition}\label{propC1}
\begin{itemize}
\item[(a)]The set $\cal P$ of all strategic measures, for a fixed initial state $x_0\in{\bf X}$, is a measurable and convex subset of ${\cal P}(\Omega)$. (Recall the notations introduced below (\ref{ea8}).)

\item[(b)]Suppose conditions (a) and (b) in Definition \ref{da12} are satisfied.  Then the space $\cal P$, endowed with the weak topology, is compact.
    \end{itemize}
\end{proposition}
\par\noindent\textit{Proof.}  For the first statement, see Theorem 8 of \cite{Piunovskiy:1997} and Chapter 5,\S5 of \cite{Dynkin:1979}. For the second statement, see e.g., \cite{Schal:1975}.
$~$\hfill$\Box$

Unless stated otherwise, we always endow the space of strategic measures with the weak topology.

The next result is known as the Derman-Strauch Lemma. It asserts that the marginal distributions of each strategy can be replicated by a Markov strategy.
\begin{proposition}\label{DTMDPChapterLemma02}
For each strategy $\sigma$, there is a Markov strategy $\sigma^M=\{\sigma^M_n\}_{n=1}^\infty$ such that
\begin{eqnarray*}
\mathsf{P}_{x_0}^\sigma(X_{n-1}\in dx,A_{n}\in da)=\mathsf{P}_{x_0}^{\sigma^M}(X_{n-1}\in dx,A_{n}\in da)
\end{eqnarray*}
for each $n=1,2,\dots.$ Here  $\sigma^M_n$ is the stochastic kernel from $\textbf{X}$ to $\textbf{A}$ such that
\begin{eqnarray*}
\mathsf{P}_{x_0}^\sigma(X_{n-1}\in dx,A_{n}\in da)= \mathsf{P}_{x_0}^{\sigma}(X_{n-1}\in dx)\sigma^M_n(da|x).
\end{eqnarray*}
One can take an arbitrarily fixed version of the stochastic kernel $\sigma^M_n$.
\end{proposition}
\par\noindent\textit{Proof.}  See  Lemma 2 of \cite{Piunovskiy:1997}. $\hfill\Box$

Now it is clear that one can restrict to Markov strategies when investigating problem (\ref{ea6}).

Next we introduce occupation measures of strategies.
\begin{definition}[Occupation measures]\label{dB1}
The {occupation measure} \index{occupation measure!in MDP} $\mathsf{M}_{x_0}^\sigma$ of a strategy $\sigma$ in the MDP $\{\textbf{X}_\Delta,\textbf{A},p\}$ with the initial state $x_0\in\textbf{X}$ is defined by
\begin{eqnarray*}
\mathsf{M}_{x_0}^\sigma(\Gamma_{X}\times \Gamma_{A})&:=& \mathsf{E}_{x_0}^\sigma\left[\sum_{n=1}^\infty I\{X_{n-1}\in \Gamma_{X},~A_{n}\in \Gamma_{A}\}\right]\\
&=&\sum_{n=1}^\infty \mathsf{E}_{x_0}^\sigma\left[I\{X_{n-1}\in\Gamma_X,~A_n\in\Gamma_A\}\right]
\end{eqnarray*}
for each $\Gamma_{X}\in {\cal B}(\textbf{X})$ and $\Gamma_{A}\in {\cal B}(\textbf{A}).$
The set of all occupation measures is denoted as ${\cal D}$; ${\cal D}^f:=\{\MM^\sigma_{x_0}:~\MM^\sigma_{x_0}({\bf X}\times{\bf A})<\infty\}$ is the set of all finite occupation measures on ${\cal B}({\bf X}\times{\bf A})$.
\end{definition}

Now  for all $j=0,1,\dots,J$,
\begin{eqnarray*}
\mathsf{E}^\sigma_{x_0}\left[\sum_{n=0}^\infty r_j(X_n,A_{n+1})\right]
=\int_{\textbf{X}\times\textbf{A}}r_j(x,a) \mathsf{M}^\sigma_{x_0}(dx\times da).
\end{eqnarray*}
Accordingly, one can reformulate problem (\ref{ea6}) as follows:
\begin{eqnarray}
 \mbox{Minimize over $\cal D:$ }~
 R_0(\MM) &:=& \int_{{\bf X}\times{\bf A}}r_0(x,a)\MM(dx\times da) \label{eB7}\\
\mbox{subject to }~R_j(\MM) &:=& \int_{{\bf X}\times{\bf A}}r_j(x,a)\MM(dx\times da)\le d_j~~~j=1,2,\ldots, J.\nonumber
\end{eqnarray}

\begin{proposition}\label{prB31}
The set of all occupation measures ${\cal D}$ with the initial state $x_0$ is a convex set in the cone of $[0,\infty]$-valued measures on ${\cal B}({\bf X}\times{\bf A})$. The set ${\cal D}^f$ of finite occupation measures is a convex subset of the linear space of finite signed measures on ${\cal B}({\bf X}\times{\bf A})$. It is a (convex) face of $\cal D$.
\end{proposition}
\par\noindent\textit{Proof.}  It follows from Proposition \ref{propC1}.  $\hfill\Box$

The next two results provide some relations satisfied by occupation measures of a strategy (respectively, a stationary strategy).
\begin{proposition}\label{DTMDPChapterThm16}
The  occupation measure $\mathsf{M}_{x_0}^\sigma$ of a strategy $\sigma$  satisfies the following equation:
\begin{eqnarray}\label{eC13}
\mu(\Gamma\times \textbf{A})=\delta_{x_0}(\Gamma)+\int_{\textbf{X}\times\textbf{A}}p(\Gamma|y,a)\mu(dy\times da),~\forall~\Gamma\in{\cal B}(\textbf{X}).
\end{eqnarray}
\end{proposition}
\par\noindent\textit{Proof.} See Lemma 9.4.3 of \cite{Hernandez-Lerma:1999}. $\hfill\Box$

\begin{proposition}\label{prB11}
Suppose $\sigma^s$ is a stationary strategy. Then
\begin{equation}\label{eB111}
\mathsf{M}_{x_0}^{\sigma^s}(\Gamma_X\times \Gamma_A)=\int_{\Gamma_X} \sigma^s(\Gamma_A|x)\mathsf{M}_{x_0}^{\sigma^s}(dx\times \textbf{A}),~\Gamma_X\in{\cal B}(\textbf{X}),~\Gamma_A\in{\cal B}(\textbf{A})
\end{equation}
and $\mathsf{M}_{x_0}^{\sigma^s}(dx\times \textbf{A})$ is the (set-wise) minimal measure on ${\cal B}(\textbf{X})$ satisfying the equation
\begin{equation}\label{eB112}
\mu(\Gamma)=\delta_{x_0}(\Gamma)+\int_{\textbf{X}}\int_{\textbf{A}} p(\Gamma|y,a)\sigma^s(da|y)\mu(dy),~\Gamma\in{\cal B}(\textbf{X}).
\end{equation}
\end{proposition}
\par\noindent\textit{Proof.} See \cite[pp.563-564]{APZY}. $\hfill\Box$

As was mentioned in Section \ref{pz23bSec01},  for discounted MDPs as well as absorbing MDPs, see e.g., \cite{feiroth,Piunovskiy:1997}, the space of occupation measures was often endowed with the weak topology generated by bounded continuous functions. To deal with infinite occupation measures, it is more convenient to endow ${\cal D}$ with the final topology generated by the projection mapping from the space of strategic measures to occupation measures. See the next definition.
\begin{definition}\label{dB3}
$\rho$ is the final topology  on $\cal D$ associated with the mapping
${O}:~{\cal P}\to{\cal D}$
defined by
$$\MM(dx\times da)=\sum_{n=1}^\infty \PP(X_{n-1}\in dx,A_n\in da).$$
That is the finest topology for which the mapping $O$ is continuous. A subset $\Gamma\subseteq {\cal D}$ is open (wrt $\rho$) if and only if $O^{-1}(\Gamma)$ is open in ${\cal P}$. Recall that ${\cal P}$ was endowed with the weak topology.
\end{definition}

\begin{lemma}\label{lB7}
Consider the MDP model $\{{\bf X},{\bf A},p\}$.
\begin{itemize}
\item[(a)] Under  conditions (a) and (b) of Definition \ref{da12}  the topological space $({\cal D},\rho)$ is compact.
\item[(b)] For each non-negative lower semicontinuous function $r(\cdot,\cdot):~{\bf X}\times{\bf A}\to\bar\RR^0_+$, the mapping $R(\cdot):~{\cal D}\to\bar\RR^0_+$  defined by
$$R(\MM):=\int_{{\bf X}\times{\bf A}} r(x,a)\MM(dx\times da)$$
is lower semicontinuous.
\end{itemize}
\end{lemma}
The proofs of all the lemmas and theorems can be found in the appendix.

\begin{corollary}\label{corB3}
If the MDP $\{{\bf X},{\bf A},p,\{r_j\}_{j=0}^J\}$ is semicontinuous, then the constrained problem
(\ref{ea6}) has an optimal solution, provided that there exists a feasible solution.
\end{corollary}
\par\noindent\textit{Proof.}
Since the equivalent  problems (\ref{ea6}) and (\ref{eB7}) have feasible solutions, the space $({\cal D},\rho)$ is compact, and the functions $R_j(\cdot)$ are lower semicontinuous, the set
$$\{\MM\in{\cal D}:~ R_j(\MM)\le d_j,~~j=1,2,\ldots, J\}$$
is nonempty and compact. Thus, the lower semicontinuous function $R_0(\cdot)$ attains its minimum thereon.  $\hfill\Box$

Alternatively, the above corollary also follows from Proposition \ref{propC1}, see also \cite{Schal:1975}, but its proof was given here in the hope of improving readability.

\section{Extreme finite occupation measures}\label{s4}
In this section we present our first main result concerning the characterization of extreme finite occupation measures. We emphasize that this result does not require any extra conditions on the MDP model; in particular, the MDP does not need to be semicontinuous.
\begin{definition}[Induced strategy]\label{dB7}
For $\MM\in{\cal D}^f$, the stationary strategy $\sigma^s$, coming form the decomposition
$$\MM(dx\times da)=\sigma^s(da|x)\MM(dx\times{\bf A})$$
on ${\cal B}({\bf X}\times{\bf A})$, is called {induced} (by $\MM$).
Here one can take an arbitrarily fixed version of the stochastic kernel $\sigma^s$, as the following lemma is valid.
\end{definition}

The next result asserts that any finite occupation measure is generated by a stationary strategy.
\begin{lemma}\label{lB8}
Suppose $\MM\in{\cal D}^f$ and $\sigma^s$ is the stationary strategy induced by $\MM$. (One can take an arbitrary version of the stochastic kernel $\sigma^s$.) Then $\MM=\MM^{\sigma^s}_{x_0}$.
\end{lemma}

The next result plays an important role in Step 1 in the proof of Theorem \ref{lB33}.
\begin{lemma}\label{lB9}
Let a stationary strategy $\sigma^s$ be such that $\MM^{\sigma^s}_{x_0}\in{\cal D}^f$. (E.g., $\sigma^s$ is the strategy,  induced by $\MM\in{\cal D}^f$.)
Then the following assertions hold.
\begin{itemize}
\item[(a)]
$$\MM^{\sigma^s}_x({\bf X}\times{\bf A})=\EE^{\sigma^s}_x\left[\sum_{n=1}^\infty I\{X_{n-1}\in{\bf X}\}\right]<\infty$$
for $\MM^{\sigma^s}_{x_0}(dx\times{\bf A})$-almost all $x\in{\bf X}$.
\item[(b)] For a bounded $\RR$-valued function $f(\cdot)$ on $\bf X$ with $f(\Delta)=0$, the function
$$v(x):=\EE^{\sigma^s}_{x}\left[\sum_{n=1}^\infty f(X_{n-1})\right]:=\EE^{\sigma^s}_{x}\left[\sum_{n=1}^\infty f^+(X_{n-1})\right]-\EE^{\sigma^s}_{x}\left[\sum_{n=1}^\infty f^-(X_{n-1})\right],~~~~~x\in{\bf X}$$
is measurable and with finite values for $\MM^{\sigma^s}_{x_0}(dx\times{\bf A})$-almost all $x\in{\bf X}$. Here the convention of $\infty-\infty:=\infty$ is in use.
\item[(c)] The function $v(\cdot)$ in (b) satisfies equation
\begin{equation}\label{eB36}
v(x)=f(x)+\int_{\bf A}\int_{\bf X} v(y)p(dy|x,a)\sigma^s(da|x)~~~~~\MM^{\sigma^s}_{x_0}(dx\times{\bf A})\mbox{-a.s.}
\end{equation}
If a measurable bounded function $w(\cdot):~{\bf X}\to\RR$ satisfies equation (\ref{eB36}), then $w(x)=v(x)$ for $\MM^{\sigma^s}_{x_0}(dx\times{\bf A})$-almost all $x\in{\bf X}$.
\end{itemize}
\end{lemma}

We note that the function $v(\cdot)$ in parts (b,c) of the previous lemma may be not finite everywhere, even though the function $f(\cdot)$ was bounded.

\begin{theorem}\label{lB33}
An occupation measure $\MM\in{\cal D}^f$  is extreme in ${\cal D}^f$ if and only if $\MM=\MM^\varphi_{x_0}$ for some deterministic stationary strategy $\varphi$.
\end{theorem}

\section{Form of the optimal control strategy}\label{s5}

In this section, we present our second main result, concerning the existence of an optimal $(J+1)$-mixed strategy to the constrained MDP problem. For this, we will impose further conditions, which, in particular, guarantee that strategies whose occupation measures are infinite are not optimal or feasible, see Theorem \ref{tB2}.

\begin{definition}\label{dB6}
According to Propositions \ref{propC1} and  \ref{prB31}, if $\sigma^1,\sigma^2,\ldots,\sigma^{\bf L}$ is a finite collection of strategies, then, for a set $\alpha_1,\alpha_2,\ldots,\alpha_{\bf L}$ of non-negative numbers with  $\sum_{l=1}^{\bf L} \alpha_l=1$,  $\sum_{l=1}^{\bf L} \alpha_l \PP_{x_0}^{\sigma^l}$ is a strategic measure and $\sum_{l=1}^{\bf L} \alpha_l \MM_{x_0}^{\sigma^l}$ is an occupation measure for some strategy $\sigma$ which is called a {mixture} of  strategies $\sigma^1,\sigma^2,\ldots,\sigma^{\bf L}$.
\end{definition}

\begin{theorem}\label{tB2}
Suppose the MDP $\{{\bf X}_\Delta,{\bf A},p,\{r_j\}_{j=0}^J\}$ with  initial state $x_0\in{\bf X}$ is  semicontinuous, and there exists a feasible strategy $\sigma$ with a finite value. Furthermore, assume that, for each strategy $\sigma$ such that $\MM^\sigma_{x_0}\notin{\cal D}^f$, there is some
$\tilde j\in\{0,1,\ldots, J\}$, possibly depending on $\sigma$, satisfying
$\int_{{\bf X}\times{\bf A}} r_{\tilde j}(x,a)\MM^\sigma_{x_0}(dx\times da)=\infty$.

Then there exists an optimal strategy in  problem (\ref{ea6}) in the form of a mixture of $J+1$ deterministic stationary strategies.
\end{theorem}

It is well known that, if the MDP is semicontinuous and the cost functions $r(\cdot,\cdot)$  are non-negative, then there exists an  optimal solution to the unconstrained problem (\ref{ea6}) (i.e., with $J=0$), which is deterministic stationary: see Corollary 9.17.2 of \cite{Bertsekas:1978} or Theorems 15.2 and 16.2 of \cite{Schal:1975a}. Therefore, we will assume that $J\ge 1$.

If there are feasible strategies in problem  (\ref{ea6}), but for all of them  $R_0(\MM^\sigma_{x_0})=+\infty$, then all feasible strategies are equally optimal. In this case, the only problem is to find a feasible strategy. To do so, we choose an arbitrary positive index, e.g., $j=1$, and investigate the problem
\begin{eqnarray*}
&\mbox{Minimize over all strategies $\sigma$:~}& \mathsf{E}_{x_0}^\sigma\left[\sum_{n=0}^\infty r_1(X_n,A_{n+1})\right]\\
&\mbox{subject to}& \mathsf{E}_{x_0}^\sigma\left[\sum_{n=0}^\infty r_j(X_n,A_{n+1})\right]\le d_j,~~~~j=2,3,\ldots, J.
\end{eqnarray*}
Clearly, after re-enumerating the indices $j$, we obtain the standard problem (\ref{ea6}) with the reduced number of constraints (or just the unconstrained problem  in case $J$ was equal to $1$). In such situations there is no need to require that the cost function $r_0(\cdot,\cdot)$ exhibits any further properties  (semicontinuity etc) except for measurability. After solving the modified problem, we obtain the desired  feasible strategy.
Clearly, the modified problem has a feasible strategy with a finite value (because the original problem had a feasible strategy).
 If all the other requirements of  Theorem \ref{tB2} are satisfied for the modified problem, then
 Theorem \ref{tB2} remains valid for it.  As the result, in such a case, there exists an optimal strategy in the original  problem (\ref{ea6}) in the form of a mixture of $J$ deterministic stationary strategies.

Let us consider the special case of optimal stopping like in \cite{dufour10}: the action space is ${\bf A}_\Delta:={\bf A}\cup\{\Delta\}$, where the isolated action $\Delta$ means stopping the process: for all $x\in{\bf X}$, $p(\{\Delta\}|x,\Delta)=1$ and $p({\bf X}|x,a)=1$ for all $a\in{\bf A}$. If this MDP $\{{\bf X}_\Delta,{\bf A}_\Delta,p,\{r_j\}_{j=0}^J\}$ is semicontinuous, there exists a feasible strategy with a finite value, and, for some $\tilde j\in\{0,1,\ldots,J\}$, $r_{\tilde j}(a,x)\ge\delta>0$ for all $x\in{\bf X}$, $a\in{\bf A}$, then all the conditions of Theorem \ref{tB2} are satisfied. $\MM^\sigma_{x_0}\notin{\cal D}^f$ means that the process is never stopped, hence $\int_{{\bf X}\times{\bf A}} r_{\tilde j}(x,a)\MM^\sigma_{x_0}(dx\times da)=\infty$. According to the above paragraph, one can omit the requirement that the feasible strategy has a finite value.

\section{Conclusion}\label{s6}

The main results of the current work are Theorems \ref{lB33} and \ref{tB2}, where we prove that every extreme finite occupation measure is generated by a deterministic stationary strategy, and, under mild conditions, show that the solution to the constrained problem is given by a finite mixture of such strategies. All the similar statements in \cite{altman,Borkar:1988,borkar91,feiroth,FeiJN,gonz,Piunovskiy:1997}, where the discounted or absorbing models were studied, follow from Theorems  \ref{lB33} and \ref{tB2}.

\section{Appendix}

\par\noindent\textit{Proof of Lemma \ref{lB7}.}
Some of the enlisted statements were presented in \cite[Lemma 4.1]{Dufour:2012}.

(a) The mapping $O$ is continuous, since ${\cal D}$ is endowed with the final topology $\rho.$ Thus, ${\cal D}={O}({\cal  P})$ is compact as the continuous image of the compact ${\cal P}$, see \cite[Chapter I, \S5, Lemma 7]{dunford}.

(b) According to Lemma 7.14(a) of \cite{Bertsekas:1978} , $r(\cdot,\cdot)=\lim_{i\to\infty} r_i(\cdot,\cdot)$, where $r_i(\cdot,\cdot)$ are point-wise increasing bounded continuous functions on ${\bf X}\times{\bf A}$. For each $i=1,2,\ldots$ the mapping
$$\PP^\sigma_{x_0}\to\int_{{\bf X}\times{\bf A}}  r_i(x,a)\PP^\sigma_{x_0}(({\bf X}\times{\bf A})^t\times dx\times da\times({\bf X}\times{\bf A})^\infty)$$
is continuous for each $t=0,1,\ldots$ because $\tau$ is the  weak topology in $\cal P$. Therefore, the mapping
\begin{eqnarray*}
\PP^\sigma_{x_0} &\to & \sum_{t=0}^n \int_{{\bf X}\times{\bf A}}  r(x,a)\PP^\sigma_{x_0}(({\bf X}\times{\bf A})^t\times dx\times da\times({\bf X}\times{\bf A})^\infty)\\
&=& \lim_{i\to\infty}\sum_{t=0}^n \int_{{\bf X}\times{\bf A}}  r_i(x,a)\PP^\sigma_{x_0}(({\bf X}\times{\bf A})^t\times dx\times da\times({\bf X}\times{\bf A})^\infty)
\end{eqnarray*}
is non-negative and lower semicontinuous again due to Lemma 7.14(a) of \cite{Bertsekas:1978}. The  monotone convergence theorem was in use here. Since $r(\cdot,\cdot)\ge 0$, Lemma B.1.1 and Proposition B.1.17 of \cite{APZY}
 imply that the mapping
\begin{eqnarray*}
\PP^\sigma_{x_0} &\to & \sum_{t=0}^\infty \int_{{\bf X}\times{\bf A}}  r(x,a)\PP^\sigma_{x_0}(({\bf X}\times{\bf A})^t\times dx\times da\times({\bf X}\times{\bf A})^\infty)\\
&=& \sup_{n=1,2,\ldots} \sum_{t=0}^n \int_{{\bf X}\times{\bf A}}  r(x,a)\PP^\sigma_{x_0}(({\bf X}\times{\bf A})^t\times dx\times da\times({\bf X}\times{\bf A})^\infty)\\
&=&\int_{{\bf X}\times{\bf A}} r(x,a)\MM^\sigma_{x_0}(dx\times da)=R(\MM^\sigma_{x_0})=R({O}(\PP^\sigma_{x_0}))
\end{eqnarray*}
is lower semicontinuous. Now, for an arbitrarily fixed $c\in\RR$
\begin{eqnarray*}
&&{O}^{-1}\left(\left\{\MM\in{\cal D}:~ R(\MM)=\int_{{\bf X}\times{\bf A}} r(x,a)\MM(dx\times da)>c\right\}\right)\\
&=&\left\{\PP\in{\cal P}:
\sum_{t=0}^\infty \int\limits_{{\bf X}\times{\bf A}}  r(x,a)\PP^\sigma_{x_0}(({\bf X}\times{\bf A})^t\times dx\times da\times({\bf X}\times{\bf A})^\infty)>c\right\}.
\end{eqnarray*}
The set  on the right-hand side is open in the topology $\tau$. Hence the set $\{\MM\in{\cal D}:~ R(M)>c\}$ is open in the topology $\rho$.

The proof is completed. $\hfill\Box$

\par\noindent\textit{Proof of Lemma \ref{lB8}.}
The both measures $\MM(dx\times{\bf A})$ and $\MM^{\sigma^s}_{x_0}(dx\times{\bf A})$ are finite and satisfy equation
\begin{equation}\label{eB33}
\mu(\Gamma_X)=\delta_{x_0}(\Gamma_X)+\int_{\bf X}\int_{\bf A} p(\Gamma_X|y,a)\sigma^s(da|y)\mu(dy),~~\forall\Gamma_X\in{\cal B}({\bf X}):
\end{equation}
see Proposition \ref{DTMDPChapterThm16} and  Proposition \ref{prB11}.

Let us show that the measure $\MM(dx\times{\bf A})$ is absolutely continuous wrt $\MM^{\sigma^s}_{x_0}(dx\times{\bf A})$ on ${\cal B}({\bf X})$.

Suppose for contradiction that $\MM(\Gamma\times{\bf A})>0$ and $\MM^{\sigma^s}_{x_0}(\Gamma\times{\bf A})=0$ for some $\Gamma\in{\cal B}({\bf X})$. Denote $\Gamma_0:=\Gamma$ and, for $n=0,1,\ldots$, put
\begin{eqnarray*}
\tilde\Gamma_{n+1} &:=& \left\{ y\in {\bf X}:~\int_{\bf A} p(\Gamma_n|y,a)\sigma^s(da|y)>0\right\};\\
\Gamma_{n+1} &:=&\Gamma_n\cup \tilde\Gamma_{n+1}.
\end{eqnarray*}
Intuitively, $\Gamma_{n+1}$ is the set of states, starting from which, the state process under $\sigma^s$  visits $\Gamma$ with positive probability within $n+1$ steps.
We will prove by induction that, for all $n=0,1,\ldots$,
\begin{eqnarray*}
&&\MM(\Gamma_n\times{\bf A})>0;\\
&&\int_{{\bf X}\setminus\Gamma_{n+1}} \int_{\bf A} p(\Gamma_n|y,a)\sigma^s(da|y)\MM(dy\times{\bf A})=0;\\
&&\MM^{\sigma^s}_{x_0}(\Gamma_n\times{\bf A})=0.
\end{eqnarray*}

When $n=0$, these assertions obviously hold because $\int_{\bf A} p(\Gamma_0|y,a)\sigma^s(da|y)=0$ for all $y\in{\bf X}\setminus \tilde\Gamma_1$. Suppose they hold for some $n\ge 0$ and consider the case of $n+1$.

Since $\Gamma_{n+1}\supseteq\Gamma_n$, $\MM(\Gamma_{n+1}\times{\bf A})>0$. Suppose $\MM^{\sigma^s}_{x_0}(\tilde\Gamma_{n+1}\times{\bf A})>0$. Then, by (\ref{eB33}),
$$\MM^{\sigma^s}_{x_0}(\Gamma_n\times{\bf A})\ge\int_{\tilde\Gamma_{n+1}}\int_{\bf A} p(\Gamma_n|y,a)\sigma^s(da|y)\MM^{\sigma^s}_{x_0}(dy\times{\bf A})>0,$$
which contradicts the inductive supposition. Thus, $\MM^{\sigma^s}_{x_0}(\tilde\Gamma_{n+1}\times{\bf A})=\MM^{\sigma^s}_{x_0}(\Gamma_{n+1}\times{\bf A})=0$. Finally,
$$\int_{{\bf X}\setminus\Gamma_{n+2}} \int_{\bf A} p(\Gamma_{n+1}|y,a)\sigma^s(da|y)\MM(dy\times{\bf A})=0$$
because $\int_{\bf A} p(\Gamma_{n+1}|y,a)\sigma^s(da|y)=0$ for all $y\in{\bf X}\setminus\tilde\Gamma_{n+2}\supseteq {\bf X}\setminus\Gamma_{n+2}$.

Therefore, for the increasing sequence $\{\Gamma_n\}_{n=0}^\infty$, after we denote $\hat\Gamma:=\bigcup_{n=0}^\infty \Gamma_n$, we have $\MM(\hat\Gamma\times{\bf A})>0$ and, by the monotone convergence theorem,
\begin{eqnarray}
&&\int_{({\bf X}\setminus\hat\Gamma)\times{\bf A}} p(\hat\Gamma|y,a)\MM(dy\times da)=\int_{{\bf X}\setminus\hat\Gamma}\int_{\bf A} p(\hat\Gamma|y,a)\sigma^s(da|y)\MM(dy\times {\bf A}) \nonumber\\
&=& \lim_{n\to\infty}  \int_{{\bf X}\setminus\hat\Gamma}\int_{\bf A} p(\Gamma_n|y,a)\sigma^s(da|y)\MM(dy\times {\bf A}) \nonumber\\
&\le & \lim_{n\to\infty}  \int_{{\bf X}\setminus\Gamma_{n+1}}\int_{\bf A} p(\Gamma_n|y,a)\sigma^s(da|y)\MM(dy\times {\bf A}) =0. \label{eB34}
\end{eqnarray}
Note also that $x_0\notin\hat\Gamma$ because $\MM^{\sigma^s}_{x_0}(\hat\Gamma\times{\bf A})=\lim_{n\rightarrow \infty} \MM^{\sigma^s}_{x_0}(\Gamma_n\times{\bf A})=0$.

Recall that $\MM\in{\cal D}^f.$ According to Proposition \ref{DTMDPChapterLemma02}, $\MM=\MM^{\sigma^M}_{x_0}$ for some Markov strategy $\sigma^M$. Since $\MM(\hat\Gamma\times{\bf A})>0$ and $x_0\notin\hat\Gamma$, there exists the minimal $n>0$ such that $\PP^{\sigma^M}_{x_0}(X_n\in\hat\Gamma)>0$, for which we have the following equalities:
\begin{eqnarray*}
0<\PP^{\sigma^M}_{x_0}(X_n\in\hat\Gamma) &=& \PP^{\sigma^M}_{x_0}(\PP^{\sigma^M}_{x_0}(X_n\in\hat\Gamma|X_{n-1}))\\
&=& \int_{\bf X}\int_{\bf A} p(\hat\Gamma|y,a)\sigma^M_n(da|y)\PP^{\sigma^M}_{x_0}(X_{n-1}\in dy)\\
&=& \int_{{\bf X}\setminus\hat\Gamma}\int_{\bf A} p(\hat\Gamma|y,a)\sigma^M_n(da|y)\PP^{\sigma^M}_{x_0}(X_{n-1}\in dy)\\
&=& \int_{({\bf X}\setminus\hat\Gamma)\times{\bf A}} p(\hat\Gamma|y,a)\PP^{\sigma^M}_{x_0}(X_{n-1}\in dy,A_n\in da),
\end{eqnarray*}
where the third equality holds by the definition of the integer $n$.
Hence,
$$\int_{({\bf X}\setminus\hat\Gamma)\times{\bf A}} p(\hat\Gamma|y,a)\MM^{\sigma^M}_{x_0}(dy\times da)=\int_{({\bf X}\setminus\hat\Gamma)\times{\bf A}} p(\hat\Gamma|y,a)\MM(dy\times da)>0,$$
which contradicts (\ref{eB34}). We have proved that $\MM(dx\times{\bf A})\ll \MM^{\sigma^s}_{x_0}(dx\times{\bf A})$.

Let us define the following substochastc kernels on ${\cal B}({\bf X})$ given $x\in{\bf X}$:
\begin{eqnarray*}
\mathsf P^0(\Gamma|x) &:= & \delta_x(\Gamma);\\
\mathsf P^{n+1}(\Gamma|x) &:=& \int_{\bf X}\int_{\bf A} p(\Gamma|y,a)\sigma^s(da|y)\mathsf P^n(dy|x),~~n=0,1,\ldots\\
&& \Gamma\in{\cal B}({\bf X}).
\end{eqnarray*}
Then $\mathsf P^i(\Gamma|x_0)=\PP^{\sigma^s}_{x_0}(X_i\in\Gamma),~i=0,1,2,\ldots. $
Now, for any finite measure $\mu$ on ${\cal B}({\bf X})$, satisfying equation  (\ref{eB33}), we have the following iterations of this equation:
\begin{eqnarray}
\mu(\Gamma) &=& \delta_{x_0}(\Gamma)+\int_{\bf X}\int_{\bf A} p(\Gamma|x_0)\sigma^s(da|x_0)\nonumber\\
&& +\int_{\bf X}\int_{\bf A} p(\Gamma|y,a)\sigma^s(da|y)\left(\int_{\bf X}\int_{\bf A} p(dy|x,a)\sigma^s(da|x)\mu(dx)\right)\nonumber\\
&=& \mathsf P^0(\Gamma|x_0)+ \mathsf P^1(\Gamma|x_0)+\int_{\bf X} \mathsf P^2(\Gamma|x)\mu(dx)\nonumber\\
&=& \mathsf P^0(\Gamma|x_0)+ \mathsf P^1(\Gamma|x_0)+\mathsf P^2(\Gamma|x_0)+\int_{\bf X} \mathsf P^3(\Gamma|x)\mu(dx)\nonumber\\
&=& \ldots  =\EE^{\sigma^s}_{x_0}\left[\sum_{i=1}^n I\{x_{i-1}\in\Gamma\}\right]+\int_{\bf X} \mathsf P^n(\Gamma|x)\mu(dx),\label{eB35}\\
&&~~~~~~~~n=1,2,\ldots.\nonumber
\end{eqnarray}
Here the Fubini Theorem was in use, and the last equality holds because
$$\mathsf P^i(\Gamma|x_0)=\PP^{\sigma^s}_{x_0}(X_i\in\Gamma),~~~i=0,1,2,\ldots.$$

Since $p({\bf X}|y,a)\le 1$, for each $x\in{\bf X}$ the sequence $\{\mathsf P^i({\bf X}|x)\}_{i=0}^\infty$ is monotonically non-increasing, so that there exists the limit $\mathsf P^\infty({\bf X}|x):=\lim_{i\to\infty}\mathsf P^i({\bf X}|x)$, and the function $\mathsf P^\infty({\bf X}|\cdot):~{\bf X}\to[0,1]$ is obviously measurable. By the dominated convergence theorem,
$$\lim_{n\to\infty} \int_{\bf X} \mathsf P^n({\bf X}|x)\mu(dx)=\int_{\bf X}\mathsf P^\infty({\bf X}|x)\mu(dx).$$
Therefore, if we substitute $\MM^{\sigma^s}_{x_0}(dx\times{\bf A})$ for $\mu(dx)$ in (\ref{eB35}), we obtain
$$\MM^{\sigma^s}_{x_0}({\bf X}\times{\bf A})=\lim_{n\to\infty} \EE^{\sigma^s}_{x_0}\left[\sum_{i=1}^n I\{X_{i-1}\in{\bf X}\}\right] + \int_{\bf X} \mathsf P^\infty({\bf X}|x)\MM^{\sigma^s}_{x_0}(dx\times{\bf A}),$$
leading to the equation $\int_{\bf X} \mathsf P^\infty({\bf X}|x)\MM^{\sigma^s}_{x_0}(dx\times{\bf A})=0$, because $\MM^{\sigma^s}_{x_0}({\bf X}\times{\bf A})$\linebreak$=\EE^{\sigma^s}_{x_0}\left[\sum_{i=1}^\infty I\{X_{i-1}\in{\bf X}\}\right]<\infty$:
the both measures $\MM^{\sigma^s}_{x_0}(dx\times{\bf A})$ and $\MM(dx\times{\bf A})$ satisfy equation (\ref{eB112}), and $\MM^{\sigma^s}_{x_0}(dx\times{\bf A})\le\MM(dx\times{\bf A})$ by Proposition \ref{prB11}. Recall that $\MM\in{\cal D}^f$.
Since $\mathsf P^\infty({\bf X}|x)\ge 0$, we conclude that $\mathsf P^\infty({\bf X}|x)=0$ $\MM^{\sigma^s}_{x_0}(dx\times{\bf A})$-a.s. and $\mathsf P^\infty({\bf X}|x)=0$ $\MM(dx\times{\bf A})$-a.s. because $\MM(dx\times{\bf A})\ll \MM^{\sigma^s}_{x_0}(dx\times{\bf A})$. Hence, for each $\Gamma\in{\cal B}({\bf X})$,
\begin{eqnarray*}
0 &\le& \limsup_{n\to\infty} \int_{\bf X} \mathsf P^n(\Gamma|x)\MM(dx\times{\bf A})\le \lim_{n\to\infty} \int_{\bf X} \mathsf P^n({\bf X}|x)\MM(dx\times{\bf A})\\
&=& \int_{\bf X} \mathsf P^\infty({\bf X}|x)\MM(dx\times{\bf A})=0,
\end{eqnarray*}
and thus $\lim_{n\to\infty} \int_{\bf X}\mathsf P^n(\Gamma|x)\MM(dx\times{\bf A})=0$. After we substitute $\MM(dx\times{\bf A})$ for $\mu(dx)$ in (\ref{eB35}), we obtain
\begin{eqnarray*}
\MM(\Gamma\times{\bf A}) &=& \lim_{n\to\infty} \EE^{\sigma^s}_{x_0}\left[\sum_{i=1}^n I\{X_{i-1}\in{\Gamma}\}\right]\\
&& +\lim_{n\to\infty} \int_{\bf X} \mathsf P^n({\Gamma}|x)\MM(dx\times{\bf A}) = \MM^{\sigma^s}_{x_0}(\Gamma\times{\bf A}).
\end{eqnarray*}
Finally,
\begin{eqnarray*}
\MM(\Gamma_X\times\Gamma_A) &=& \int_{\bf X} \sigma^s(\Gamma_A|x)\MM(dx\times{\bf A})\\
&=&  \int_{\bf X} \sigma^s(\Gamma_A|x)\MM^{\sigma^s}_{x_0}(dx\times{\bf A})=\MM^{\sigma^s}_{x_0}(\Gamma_X\times\Gamma_A),\\
&&~~~~~~~~\forall \Gamma_X\in{\cal B}({\bf X}),~\Gamma_A\in{\cal B}({\bf A}).
\end{eqnarray*}
$\hfill\Box$

\par\noindent\textit{Proof of Lemma \ref{lB9}.} Note that, if a statement  $S(X_m)$ is valid $\PP^{\sigma^s}_{x_0}$-a.s. for all $m=0,1,2,\ldots$, then the statement $S(x)$ is  valid for $\MM^{\sigma^s}_{x_0}(dx\times{\bf A})$-almost all $x\in{\bf X}$ and vice versa.

(a)  If the formulated statement does not hold, then there is a set $\Gamma\in{\cal B}({\bf X})$ such that, for some $m\ge 0$, $\PP^{\sigma^s}_{x_0}(X_m\in\Gamma)>0$ and
$$\EE^{\sigma^s}_{x}\left[\sum_{n=1}^\infty I\{X_{n-1}\in{\bf X}\}\right]=\infty~~~\forall x\in\Gamma.$$
Now
\begin{eqnarray*}
\MM^{\sigma^s}_{x_0}({\bf X}\times{\bf A})&\ge & \EE^{\sigma^s}_{x_0}\left[\sum_{n=m+1}^\infty I\{X_{n-1}\in{\bf X}\}\right]\\
&=&\EE^{\sigma^s}_{x_0}\left[\EE^{\sigma^s}_{x_0}\left[\left.\sum_{n=m+1}^\infty I\{X_{n-1}\in{\bf X}\}\right|X_m\right]\right]\\
&\ge & \int_\Gamma \EE^{\sigma^s}_x \left[\sum_{n=1}^\infty I\{X_{n-1}\in{\bf X}\}\right] \PP^{\sigma^s}_{x_0}(X_m\in dx)=+\infty.
\end{eqnarray*}
Here
$$\EE^{\sigma^s}_{x_0}\left[\left.\sum_{n=m+1}^\infty I\{X_{n-1}\in{\bf X}\}\right|X_m\right]=\EE^{\sigma^s}_{X_m}\left[\sum_{n=1}^\infty I\{X_{n-1}\in{\bf X}\}\right]$$
because the controlled process $\{X_n\}_{n=0}^\infty$, under the strategy $\sigma^s$, is Markov and time-homogeneous.

The obtained contradiction with the assumption $\MM^{\sigma^s}_{x_0}\in{\cal D}^f$ proves the statement.

(b) Let $\{{\cal F}_n\}_{n=0}^\infty$ be the natural filtration ${\cal F}_n:=\sigma\{X_0,X_1,\ldots,X_n\}$ of the Markov time-homogeneous process $\{X_n\}_{n=0}^\infty$ under the control strategy $\sigma^s$ and with the initial state $x_0\in{\bf X}$.
For the positive and negative parts of $f(\cdot)$, we have the following relations for each fixed $m\ge 0$:
\begin{eqnarray*}
0 &\le & \EE^{\sigma^s}_{x_0}\left[\left.\sum_{n=m+1}^\infty f^{\pm}(X_{n-1})\right|{\cal F}_m\right]\\
&=& \EE^{\sigma^s}_{X_m}\left[\sum_{n=1}^\infty f^{\pm}(X_{n-1})\right]<\infty~~~~~~~\PP^{\sigma^s}_{x_0}\mbox{-a.s.}
\end{eqnarray*}
The very last inequality holds by (a) and the boundedness of the function $f(\cdot)$. Therefore, the function $v(\cdot)$ is with finite values for $\MM^{\sigma^s}_{x_0} (dx\times{\bf A})$-almost all $x\in{\bf X}$.

The measurability of $v(\cdot)$ follows from Proposition \ref{propC1}(b). Statement (b) is proved.

(c)
Let $\{{\cal F}_n\}_{n=0}^\infty$ be the natural filtration ${\cal F}_n:=\sigma\{X_0,X_1,\ldots,X_n\}$ of the Markov time-homogeneous process $\{X_n\}_{n=0}^\infty$ under the control strategy $\sigma^s$ and with the initial state $x\in{\bf X}$.
According to (b),
\begin{eqnarray*}
v(x)&=&f(x)+\EE^{\sigma^s}_x\left[\EE^{\sigma^s}_x\left[\left.\sum_{n=2}^\infty f(X_{n-1})\right|{\cal F}_1\right]\right]\\
&=& f(x)+\EE^{\sigma^s}_x[v(X_1)]\\
&=&f(x)+\int_{\bf A}\int_{\bf X} v(y)p(dy|x,a)\sigma^s(da|x)~~~\MM^{\sigma^s}_{x_0}(dx\times{\bf A})\mbox{-a.s.}
\end{eqnarray*}
Here, all the terms are finite $\MM^{\sigma^s}_{x_0}(dx\times{\bf A})$-a.s. by (b). Equation (\ref{eB36}) is proved.

Let us fix  an arbitrary $i\in\{0,1,\ldots\}$ and the filtration $\{{\cal F}_n\}_{n=0}^\infty$ corresponding to the initial state $x_0\in{\bf X}$. For the function $w(\cdot)$, we have the following obvious equation
\begin{equation}\label{eB38}
\EE^{\sigma^s}_{x_0}[w(X_{i+1})|{\cal F}_i]=\int_{\bf A}\int_{\bf X} w(y)p(dy|X_i,a)\sigma^s(da|X_i)~~\PP^{\sigma^s}_{x_0}\mbox{-a.s.}
\end{equation}
We are going to prove by induction the following equality
\begin{equation}\label{eB39}
w(X_i)=\EE^{\sigma^s}_{x_0}\left[\left.\sum_{j=0}^k f(X_{i+j})\right|{\cal F}_i\right]+\EE^{\sigma^s}_{x_0}[w(X_{i+k+1})|{\cal F}_i]~~\PP^{\sigma^s}_{x_0}\mbox{-a.s.}
\end{equation}
for $k=0,1,\ldots$.

When $k=0$, equality (\ref{eB39}) follows from equation (\ref{eB36}) for $w(\cdot)$ and (\ref{eB38}). Suppose it holds for some $k\ge 0$. Then
\begin{eqnarray*}
w(X_i) &=& \EE^{\sigma^s}_{x_0}\left[\left. \sum_{j=0}^k f(X_{i+j})\right|{\cal F}_i\right]+ \EE^{\sigma^s}_{x_0}[f(X_{i+k+1})|{\cal F}_i]\\
&&+\EE^{\sigma^s}_{x_0}\left[\left.\int_{\bf A}\int_{\bf X} w(y)p(dy|X_{i+k+1},a)\sigma^s(da|X_{i+k+1})\right|{\cal F}_i\right]\\
&=&  \EE^{\sigma^s}_{x_0}\left[\left. \sum_{j=0}^{k+1} f(X_{i+j})\right|{\cal F}_i\right]+\EE^{\sigma^s}_{x_0}[\EE^{\sigma^s}_{x_0}[w(X_{i+k+2})|{\cal F}_{i+k+1}]|{\cal F}_i]\\
&=&  \EE^{\sigma^s}_{x_0}\left[\left. \sum_{j=0}^{k+1} f(X_{i+j})\right|{\cal F}_i\right]+\EE^{\sigma^s}_{x_0}[w(X_{i+k+2})|{\cal F}_i]~~~\PP^{\sigma^s}_{x_0}\mbox{-a.s.}
\end{eqnarray*}
Here, the first equality is by (\ref{eB36}), and the second equality is by (\ref{eB38}). Equality (\ref{eB39}) is proved.

When $k\to\infty$, since the process $\{X_n\}_{n=0}^\infty$ is Markov and time-homogeneous,
$$ \EE^{\sigma^s}_{x_0}\left[\left. \sum_{j=0}^{k} f(X_{i+j})\right|{\cal F}_i\right]= \EE^{\sigma^s}_{X_i}\left[ \sum_{n=1}^{k+1} f(X_{n-1})\right]\to v(X_i)~~\PP^{\sigma^s}_{x_0}\mbox{-a.s.}$$
by the definition of the function $v(\cdot)$. According to (a),
$$\EE^{\sigma^s}_{x_0}[I\{X_{i+k+1}\in{\bf X}\}|{\cal F}_i]=\EE^{\sigma^s}_{X_i}[I\{X_{k+1}\in{\bf X}\}]\to 0 \mbox{ as } k\to\infty ~\PP^{\sigma^s}_{x_0}\mbox{-a.s.}$$
Therefore, since the function $w(\cdot)$ is bounded, $w(X_i)=v(X_i)$ $\PP^{\sigma^s}_{x_0}$-a.s.

The proof is completed.
$\hfill\Box$

\par\noindent\textit{Proof of Theorem \ref{lB33}.}  We assume that ${\cal D}^f\ne\emptyset$ and, according to  Lemma \ref{lB8},   consider only  the occupation measures $\MM=\MM^{\sigma^s}_{x_0}\in{\cal D}^f$ coming from stationary strategies $\sigma^s$.

(a) The `if' part.  We will prove a little more general statement: if $\MM^\varphi_{x_0}\in{\cal D}^f$ is the occupation measure generated by a deterministic stationary strategy $\varphi$, then $\MM^\varphi_{x_0}$ is extreme in $\cal D$ (and certainly in ${\cal D}^f$, too).

Suppose
$\MM^\varphi_{x_0}=\alpha \MM_1+(1-\alpha)\MM_2$ with $\alpha\in(0,1)$ and $\MM_{1,2}\in{\cal D}$. Then $\MM_{1,2}\in{\cal D}^f$ because $\MM^\varphi_{x_0}\in{\cal D}^f$ and, according to Lemma \ref{lB8}, $\MM_{1,2}=\MM^{\sigma^s_{1,2}}_{x_0}$ for the induced stationary strategies $\sigma^s_{1,2}$. Therefore,
\begin{equation}\label{eB1131}
\MM=\MM^\varphi_{x_0}=\alpha\MM^{\sigma^s_1}_{x_0}+(1-\alpha)\MM^{\sigma^s_2}_{x_0}.
\end{equation}
The goal is to show that $\MM^{\sigma^s_1}_{x_0}=\MM^{\sigma^s_2}_{x_0}=\MM^\varphi_{x_0}$.

The both marginal measures $\MM^{\sigma^s_1}_{x_0}(dx\times{\bf A})$ and  $\MM^{\sigma^s_2}_{x_0}(dx\times{\bf A})$ are absolutely continuous wrt $\MM^{\varphi}_{x_0}(dx\times{\bf A})$; the Radon-Nikodym derivatives are denoted as $h_1(\cdot)$ and $h_2(\cdot)$ correspondingly. From (\ref{eB1131}) we have
\begin{equation}\label{eB1141}
\alpha h_1(x)+(1-\alpha) h_2(x)=1 \mbox{ for } \MM^{\varphi}_{x_0}(dx\times{\bf A})\mbox{-almost all } x\in{\bf X}.
\end{equation}
Now, using (\ref{eB111}), we have equalities
\begin{eqnarray*}
\MM^{\varphi}_{x_0}(dx\times da) &=& \MM^{\varphi}_{x_0}(dx\times{\bf A})\delta_{\varphi(x)}(da)\\
&=& \alpha \MM^{\varphi}_{x_0}(dx\times{\bf A})h_1(x)\sigma^s_1(da|x)\\
&&+(1-\alpha)\MM^{\varphi}_{x_0}(dx\times{\bf A})h_2(x)\sigma^s_2(da|x)\\
&=& \MM^{\varphi}_{x_0}(dx\times{\bf A})[\alpha h_1(x)\sigma^s_1(da|x)+(1-\alpha) h_2(x) \sigma^s_2(da|x)].
\end{eqnarray*}
The expression in the square brackets is the Dirac measure $\delta_{\varphi(x)}(da)$ for $\MM^{\varphi}_{x_0}(dx\times{\bf A})$-almost all  $x\in{\bf X}$. Note that any Dirac measure on ${\cal B}({\bf A})$ is extreme in ${\cal P}({\bf A})$. Therefore, using (\ref{eB1141}), we conclude that
\begin{itemize}
\item on the set ${\bf I}_0:=\{x\in{\bf X}:~\alpha h_1(x)\in(0,1)\}$, $\sigma^s_1(da|x)=\sigma^s_2(da|x)=\delta_{\varphi(x)}(da)$ for  $\MM^{\varphi}_{x_0}(dx\times{\bf A})$-almost all  $x\in{\bf I}_0$;
\item on the set ${\bf I}_1:=\{x\in{\bf X}:~\alpha h_1(x)=1\}$, $\sigma^s_1(da|x)=\delta_{\varphi(x)}(da)$ for  $\MM^{\varphi}_{x_0}(dx\times{\bf A})$-almost all  $x\in{\bf I}_1$, and the stochastic kernel $\sigma^s_2(da|x)$ may be arbitrary, but on the set ${\bf I}_1$, since $(1-\alpha)>0$, $h_2(x)=0$ for  $\MM^{\varphi}_{x_0}(dx\times{\bf A})$-almost all  $x\in{\bf I}_1$, i.e., $\MM^{\sigma^s_2}_{x_0}({\bf I}_1\times{\bf A})=0$, and the values of $\sigma^s_2(da|x)$ are of no importance for the measure $\MM^{\sigma^s_2}_{x_0}(dx\times da)$ on ${\cal B}({\bf I}_1\times{\bf A})$;
\item symmetrically, on the set ${\bf I}_2:=\{x\in{\bf X}:~(1-\alpha) h_2(x)=1\}$, $\sigma^s_2(da|x)=\delta_{\varphi(x)}(da)$ for  $\MM^{\varphi}_{x_0}(dx\times{\bf A})$-almost all  $x\in\{{\bf I}_2$ and $\MM^{\sigma^s_1}_{x_0}({\bf I}_2\times{\bf A})=0$.
\end{itemize}
Recall that the measures $\MM^{\sigma^s_1}_{x_0}(dx\times{\bf A})$ and  $\MM^{\sigma^s_2}_{x_0}(dx\times{\bf A})$ are absolutely continuous wrt $\MM^{\varphi}_{x_0}(dx\times{\bf A})$. Thus, $\sigma^s_1(da|x)=\delta_{\varphi(x)}(da)$ for $\MM^{\sigma^s_1}_{x_0}(dx\times{\bf A})$-almost all $x\in{\bf X}$ and $\sigma^s_2(da|x)=\delta_{\varphi(x)}(da)$ for $\MM^{\sigma^s_2}_{x_0}(dx\times{\bf A})$-almost all $x\in{\bf X}$.
 Hence, $\varphi$ is the stationary strategy induced by $\MM^{\sigma^s_{1,2}}_{x_0}$. and
$$\MM^{\sigma^s_1}_{x_0}=\MM^{\sigma^s_2}_{x_0}=\MM^\varphi_{x_0}.$$
by Lemma \ref{lB8}.

(b) The `only if' part. According to Lemma  \ref{lB8}, an extreme point $\MM$ in ${\cal D}^f$ satisfies the equalities  $\MM(dx\times da)=\MM^{\sigma^s}_{x_0}(dx\times da)=\sigma^s(da|x)\MM^{\sigma^s}_{x_0}(dx\times{\bf A})$ on ${\cal B}({\bf X}\times{\bf A})$ for the induced stationary strategy $\sigma^s$.

\underline{Step 1.} Suppose that
$$\sigma^s(da|x)=\alpha\sigma^s_1(da|x)+(1-\alpha) \sigma^s_2(da|x),$$
where $\alpha\in(0,1)$ and $\sigma^s_1$ and $\sigma^s_2$ are two essentially different  stochastic kernels on $\bf A$ given $\bf X$. To be precise, we assume that,
for some $\hat\Gamma^A\in{\cal B}({\bf A})$ and $\hat\Gamma^X\in{\cal B}({\bf X})$,
$$\MM^{\sigma^s}_{x_0}(\hat\Gamma^X\times{\bf A})>0 \mbox{ and } \sigma^s_2(\hat\Gamma^A|x)>\sigma^s_1(\hat\Gamma^A|x) \mbox{ for all } x\in\hat\Gamma^X.$$
The stochastic kernels $\sigma^s_{1,2}$ define  the corresponding stationary strategies, again denoted as $\sigma^s_{1,2}$.
We will show that, in this case,  the measure $\MM^{\sigma^s}_{x_0}$ is not extreme in ${\cal D}^f$.

If $\MM^{\sigma^s_1}_{x_0}(dx\times{\bf A})=\MM^{\sigma^s_2}_{x_0}(dx\times{\bf A})=\MM^{\sigma^s}_{x_0}(dx\times{\bf A})$, then, by Lemma \ref{prB11},
\begin{eqnarray*}
\MM^{\sigma^s}_{x_0}(dx\times da) &=& \MM^{\sigma^s}_{x_0}(dx\times{\bf A})\sigma^s(da|x)\\
&=& \alpha \MM^{\sigma^s_1}_{x_0}(dx\times{\bf A})\sigma^s_1(da|x)+(1-\alpha)\MM^{\sigma^s_2}_{x_0}(dx\times{\bf A})\sigma_2(da|x)\\
&=& \alpha \MM^{\sigma^s_1}_{x_0}(dx\times da)+(1-\alpha)\MM^{\sigma^s_2}_{x_0}(dx\times da).
\end{eqnarray*}
Therefore, the measure $\MM^{\sigma^s}_{x_0}$ is not extreme in ${\cal D}^f$, as $\MM^{\sigma^s_1}_{x_0}\ne \MM^{\sigma^s_2}_{x_0}$ and $\MM^{\sigma^s_1}_{x_0}, \MM^{\sigma^s_2}_{x_0}\in{\cal D}^f$. (Recall that ${\cal D}^f$ is a face of $\cal D$.)

Suppose now that $\MM^{\sigma^s_1}_{x_0}(dx\times{\bf A})\ne\MM^{\sigma^s}_{x_0}(dx\times{\bf A})$ or $\MM^{\sigma^s_2}_{x_0}(dx\times{\bf A})\ne\MM^{\sigma^s}_{x_0}(dx\times{\bf A})$ . There exists the first moment $\tau>0$ such that
$$\mbox{either } \PP^{\sigma^s_1}_{x_0}(X_\tau \in dx)\ne\PP^{\sigma^s}_{x_0}(X_\tau\in dx)~\mbox{ or }  \PP^{\sigma^s_2}_{x_0}(X_\tau \in dx)\ne\PP^{\sigma^s}_{x_0}(X_\tau\in dx).$$
Without loss of generality we assume that the first inequality holds. What actually happens is that the both inequalities hold simultaneously at the moment $\tau$: see (\ref{eB1211}).

Since
$$\PP^{\sigma^s_1}_{x_0}(X_{\tau-1}\in dx)=\PP^{\sigma^s_2}_{x_0}(X_{\tau-1}\in dx)=\PP^{\sigma^s}_{x_0}(X_{\tau-1}\in dx),$$
we have equality
\begin{equation}\label{eB391}
\PP^{\sigma^s}_{x_0}(X_\tau\in dx)=\alpha\PP^{\sigma^s_1}_{x_0}(X_{\tau}\in dx)+(1-\alpha)\PP^{\sigma^s_2}_{x_0}(X_{\tau}\in dx)
\end{equation}
because
$$\sigma^s(da|x)=\alpha\sigma^s_1(da|x)+(1-\alpha) \sigma^s_2(da|x).$$
Note, the controlled process $\{X_n\}_{n=0}^\infty$ is Markov and time-homogeneous under all strategies $\sigma^s$ and $\sigma^s_{1,2}$, with the transition probabilities
$$\int_{\bf A} p(dy|x,a)\sigma^s(da|x)~~\mbox{and }~\int_{\bf A} p(dy|x,a)\sigma^s_{1,2}(da|x)$$
correspondingly.

Let us introduce the Markov (non-stationary) strategies $\sigma^{M_1}$ and $\sigma^{M_2}$ by the formulae
$$\sigma^{M_{1,2}}_n(da|x)=I\{n\ne\tau\}\sigma^s(da|x)+I\{n=\tau\} \sigma^s_{1,2}(da|x).$$

The combination of strategic measures $\alpha\PP^{\sigma^{M_1}}_{x_0}+(1-\alpha)\PP^{\sigma^{M_2}}_{x_0}$ satisfies the key properties of the strategic measure $\PP^{\sigma^s}_{x_0}$: see (\ref{ea3}),(\ref{ea4}),(\ref{ea8}), or formula (1.7) in \cite{Piunovskiy:1997}, or \cite[\S2.2.3]{Hernandez-Lerma:1996}. Thus,
\begin{equation}\label{eB161}
\PP^{\sigma^s}_{x_0}=\alpha\PP^{\sigma^{M_1}}_{x_0}+(1-\alpha)\PP^{\sigma^{M_2}}_{x_0}\Longrightarrow \MM^{\sigma^s}_{x_0}=\alpha\MM^{\sigma^{M_1}}_{x_0}+(1-\alpha)\MM^{\sigma^{M_2}}_{x_0}
\end{equation}
and, like previously, $\MM^{\sigma^{M_1}}_{x_0},\MM^{\sigma^{M_2}}_{x_0}\in{\cal D}^f$ because $\MM^{\sigma^s}_{x_0}\in{\cal D}^f$. We aim to show that $\MM^{\sigma^{M_1}}_{x_0}(dx\times{\bf A})\ne \MM^{\sigma^{M_2}}_{x_0}(dx\times{\bf A})$  on ${\cal B}({\bf X})$, leading to the desired assertion that $\MM^{\sigma^s}_{x_0}$ is not extreme in ${\cal D}^f$.

Since the strategies $\sigma^s$, $\sigma^s_{1,2}$ and $\sigma^{M_{1,2}}$ are Markov and $\sigma^s=\alpha\sigma^s_1+(1-\alpha)\sigma^s_2$, we have the following relations:
\begin{eqnarray}
\PP^{\sigma^s}_{x_0}(X_\tau\in dx) &=&\alpha \PP^{\sigma^{M_1}}_{x_0}(X_\tau\in dx)+(1-\alpha) \PP^{\sigma^{M_2}}_{x_0}(X_\tau\in dx)\nonumber\\
&=& \alpha \PP^{\sigma^s_{1}}_{x_0}(X_\tau\in dx)+(1-\alpha) \PP^{\sigma^s_{2}}_{x_0}(X_\tau\in dx);\nonumber \\
\PP^{\sigma^{M_1}}_{x_0}(X_\tau\in dx) &= &  \PP^{\sigma^s_{1}}_{x_0}(X_\tau\in dx)\ne  \PP^{\sigma^{s}}_{x_0}(X_\tau\in dx); \nonumber\\
\PP^{\sigma^{M_2}}_{x_0}(X_\tau\in dx) & = &  \PP^{\sigma^s_{2}}_{x_0}(X_\tau\in dx) \ne  \PP^{\sigma^{M_1}}_{x_0}(X_\tau\in dx),~~ \PP^{\sigma^{s}}_{x_0}(X_\tau\in dx).\nonumber\\
&& \label{eB1211}
\end{eqnarray}
The first three lines here are according to the definitions of $\tau$ and of the strategies $\sigma^{M_{1,2}}$ (see also (\ref{eB391})), and  inequalities (\ref{eB1211}) follow from them.

Let $\Gamma\in{\cal B}({\bf X})$ be such that $\PP^{\sigma^{M_1}}_{x_0}(X_\tau\in\Gamma)\ne \PP^{\sigma^{M_2}}_{x_0}(X_\tau\in\Gamma)$.
We fix the following bounded functions on ${\bf X}_\Delta$, equal to zero on $\Delta$:
\begin{eqnarray*}
h(x) &:=& I\{x\in\Gamma\};\\
f(x) &:=& h(x)-\int_{\bf A}\int_{\bf X} h(y) p(dy|x,a)\sigma^s(da|x).
\end{eqnarray*}

According to Lemma \ref{lB9}(b,c), the function
$$v(x):=\EE^{\sigma^s}_{x} \left[\sum_{n=1}^\infty f(X_{n-1})\right],~~~~~x\in{\bf X}$$
is measurable and equals $h(x)$ for $\MM^{\sigma^s}_{x_0}(dx\times{\bf A})$-almost all $x\in{\bf X}$ because $h(\cdot)$ satisfies equation (\ref{eB36}).

Let $\{{\cal F}\}_{t=0}^\infty$ be the natural filtration of the process $\{X_n\}_{n=0}^\infty$, i.e., ${\cal F}_t:=\sigma\{X_0,X_1,\ldots, $\linebreak$X_t\}$. According to the definition of the strategies $\sigma^{M_{1,2}}$, $\sigma^{M_{1,2}}_n=\sigma^s$ for $n>\tau$, so
\begin{eqnarray*}
&&\EE^{\sigma^{M_{1,2}}}_{x_0}\left[\left.\sum_{n=\tau+1}^\infty f(X_{n-1})\right|{\cal F}_\tau\right] = \EE^{\sigma^s}_{x_0}\left[\left.\sum_{n=\tau+1}^\infty f(X_{n-1})\right|{\cal F}_\tau\right]\\
&=& \EE^{\sigma^s}_{X_\tau}\left[\sum_{n=1}^\infty f(X_{n-1})\right]=h(X_\tau)~~~\PP^{\sigma^s}_{x_0}\mbox{-a.s., \ and thus } \PP^{\sigma^{M_{1,2}}}_{x_0}\mbox{-a.s.}
\end{eqnarray*}
The second equality holds because the controlled process $\{X_n\}_{n=0}^\infty$ under the stationary strategy $\sigma^s$ is Markov and time-homogeneous. Note also that $\PP^{\sigma^{M_{1,2}}}_{x_0}\ll \PP^{\sigma^s}_{x_0}$ by (\ref{eB161}). Now, since $\sigma^{M_{1,2}}_n=\sigma^s$ for $n<\tau$,
\begin{eqnarray*}
\int_{\bf X} f(x)\MM^{\sigma^{M_{1,2}}} _{x_0}(dx\times{\bf A}) &=& \EE^{\sigma^s}_{x_0}\left[\sum_{n=1}^\tau f(X_{n-1})\right]\\
&&+\EE^{\sigma^{M_{1,2}}}_{x_0}\left[\EE^{\sigma^{M_{1,2}}}_{x_0}\left[\left. \sum_{n=\tau+1}^\infty f(X_{n-1})\right|{\cal F}_\tau\right]\right]\\
&=& \EE^{\sigma^s}_{x_0}\left[\sum_{n=1}^\tau f(X_{n-1})\right]+\EE^{\sigma^{M_{1,2}}}_{x_0}[h(X_\tau)]\\
&=& \EE^{\sigma^s}_{x_0}\left[\sum_{n=1}^\tau f(X_{n-1})\right]+\PP^{\sigma^{M_{1,2}}}_{x_0}(X_\tau\in\Gamma).
\end{eqnarray*}
As the result, by the definition of the subset $\Gamma$,
\begin{eqnarray*}
&& \int_{\bf X} f(x)\MM^{\sigma^{M_{1}}} _{x_0}(dx\times{\bf A}) \ne \int_{\bf X} f(x)\MM^{\sigma^{M_{2}}} _{x_0}(dx\times{\bf A})\\
&\Longrightarrow & \MM^{\sigma^{M_{1}}} _{x_0}(dx\times{\bf A}) \ne \MM^{\sigma^{M_{2}}} _{x_0}(dx\times{\bf A})~~\mbox{ on }~{\cal B}({\bf X}).
\end{eqnarray*}
Hence, the measure $\MM^{\sigma^s}_{x_0}$ is not extreme in ${\cal D}^f$.

The further steps in fact repeat the proof of Theorem 10 of \cite{Piunovskiy:1997}, but we provide the details for completeness.

\underline{Step 2.} We will show that,  if $\MM^{\sigma^s}_{x_0}$ is an extreme point in ${\cal D}^f$, then, for each $\Gamma^A\in{\cal B}({\bf A})$,   $\Gamma^X\in{\cal B}({\bf X})$, $\alpha\in(0,1)$, in case $\MM^{\sigma^s}_{x_0}(\Gamma^X\times{\bf A})>0$, there is $x\in\Gamma^X$ such that either $\sigma^s(\Gamma^A|x)<\alpha$ or $\sigma^s(\Gamma^A|x)>1-\alpha$.

This statement is trivial for $\alpha > 1/2$: if $\sigma^s(\Gamma^A|x)\ge\alpha>1/2$, then $1-\sigma^s(\Gamma^A|x)<1/2<\alpha$. Thus, below we assume that $\alpha\in(0,1/2]$.

The proof is  by contradiction. Namely, suppose there exist $\hat\Gamma^A\in{\cal B}({\bf A})$, $\hat\Gamma^X\in{\cal B}({\bf X})$ and $\alpha\in(0,1/2]$ such that $\MM^{\sigma^s}_{x_0}(\hat\Gamma^X\times{\bf A})>0$ and $\sigma^s(\hat\Gamma^A|x)\in[\alpha,1-\alpha]$ for all $x\in\hat\Gamma^X$. Consider the following stochastic kernels:
\begin{eqnarray*}
\sigma^s_1(\Gamma^A|x) &=& \left\{\begin{array}{ll}
\sigma^s(\Gamma^A|x), & \mbox{if } x\notin\hat\Gamma^X;\\
\frac{\sigma^s(\Gamma^A\cap(\hat\Gamma^A)^c|x)}{\sigma^s((\hat\Gamma^A)^c|x)}, & \mbox{if }x\in\hat\Gamma^X;
\end{array}\right.\\
\sigma^s_2(\Gamma^A|x) &=& \left\{\begin{array}{l}
\sigma^s(\Gamma^A|x), \hspace{10mm} \mbox{if } x\notin\hat\Gamma^X;\\
\frac{\sigma^s(\Gamma^A\cap\hat\Gamma^A|x)}{1-\alpha}+ \sigma^s(\Gamma^A\cap(\hat\Gamma^A)^c|x)\frac{ \sigma^s((\hat\Gamma^A)^c|x)-\alpha}{(1-\alpha) \sigma^s((\hat\Gamma^A)^c|x)}, \\
~\hspace{25mm} \mbox{if } x\in\hat\Gamma^X;
\end{array}\right.\\
&&\Gamma^A\in{\cal B}({\bf A}),
\end{eqnarray*}
which are well defined because $\sigma^s((\hat\Gamma^A)^c|x)\ge \alpha>0$.

Clearly, $\alpha\sigma^s_1(\Gamma^A|x)+(1-\alpha)\sigma^s_2(\Gamma^A|x)=\sigma^s(\Gamma^A|x)$ for all $\Gamma^A\in{\cal B}({\bf A})$ and $x\in{\bf X}$; $\MM^{\sigma^s}_{x_0}(\hat\Gamma^X\times{\bf A})>0$, and, for all $x\in\hat\Gamma^X$,
$$\sigma^s_2(\hat\Gamma^A|x)-\sigma^s_1(\hat\Gamma^A|x)=\frac{\sigma^s(\hat\Gamma^A|x)}{1-\alpha}-0\ge \frac{\alpha}{1-\alpha}>0.$$
Therefore, the measure $\MM^{\sigma^s}_{x_0}$ is not extreme in ${\cal D}^f$ according to the statement in Step 1.

\underline{Step 3.} We will show that,  if $\MM^{\sigma^s}_{x_0}$ is an extreme point in ${\cal D}^f$, then, for each $\Gamma^A\in{\cal B}({\bf A})$,
$$\sigma^s(\Gamma^A|x)\in\{0,1\} \mbox{ for } \MM^{\sigma^s}_{x_0}(dx\times{\bf A})\mbox{-almost all } x\in{\bf X}.$$

Let $\Gamma^A\in{\cal B}({\bf A})$ be arbitrarily fixed and introduce the sets
$$\Gamma^X_{\Gamma^A}(i):=\left\{x\in{\bf X}:~\sigma^s(\Gamma^A|x)\in\left[\left(\frac{1}{2}\right)^i,1-\left(\frac{1}{2}\right)^i\right]\right\}\in{\cal B}({\bf X}),~i=1,2,\ldots.$$
For eah $i=1,2,\ldots$, $\MM^{\sigma^s}_{x_0}(\Gamma^X_{\Gamma^A}(i) \times{\bf A})=0$ because, otherwise, for $\Gamma^A,\Gamma^X:=\Gamma^X_{\Gamma^A}(i)$, and $\alpha:=(\frac{1}{2})^i$, we would have $\MM^{\sigma^s}_{x_0}(\Gamma^X\times{\bf A})>0$ and, for all $x\in\Gamma^X$, $\sigma^s(\Gamma^A|x)\in[\alpha,1-\alpha]$, which contradicts the statement proved at Step 2.

Note that $\Gamma^X_{\Gamma^A}(i)\subset \Gamma^X_{\Gamma^A}(i+1)$ for all $i=1,2,\ldots$. Now, for
$$\Gamma^X_{\Gamma^A}:=\bigcup_{i=1}^\infty \Gamma^X_{\Gamma^A}(i)=\{x\in{\bf X}:~\sigma^s(\Gamma^A|x)\in(0,1)\},$$
we have $\MM^{\sigma^s}_{x_0}(\Gamma^X_{\Gamma^A}\times{\bf A})=\lim_{i\to\infty}\MM^{\sigma^s}_{x_0}(\Gamma^X_{\Gamma^A}(i)\times{\bf A})=0$,
and the desired statement is proved.

\underline{Step 4.} Finally, we proceed to construct the measurable mapping $\varphi:~{\bf X}\to{\bf A}$ such that
$$\sigma^s(da|x)=\delta_{\varphi(x)}(da) \mbox{ for } \MM^{\sigma^s}_{x_0}(dx\times{\bf A})\mbox{-almost all } x\in{\bf X}.$$

As a separable metrizable space, $\bf A$ has a totally bounded metrization $\kappa$ \cite[Corollary 7.6.1]{Bertsekas:1978}:
$$\forall \varepsilon>0~~\exists \{a_1,a_2,\ldots, a_n\}\subset{\bf A}:~~{\bf A}=\bigcup_{i=1}^n O(a_i,\varepsilon),$$
where $O(a_i,\varepsilon):=\{a\in{\bf A}:~\kappa(a,a_i)<\varepsilon\}$.

Let $\varepsilon_k:=\frac{1}{2^k}$ ($k=1,2,\ldots$) and let $\{a^k_1,a^k_2,\ldots, a^k_{n_k}\}$ be the corresponding $\varepsilon_k$-net in $\bf A$. Denote
$$S^k_i:=\{x\in{\bf X}:~\sigma^s(O(a^k_i,\varepsilon_k)|x)\notin\{0,1\}\}$$
and $S:=\bigcup_{k,i}S^k_i$. These sets are obviously measurable, and, for all $k=1,2,\ldots,$ $i=1,2,\ldots, n_k$, $\MM^{\sigma^s}_{x_0}(S^k_i\times{\bf A})=0$ according to the statement proved on Step 3; $\MM^{\sigma^s}_{x_0}(S\times{\bf A})=0$, too.

We are going to construct the desired mapping $\varphi$ on ${\bf X}\setminus S$ and then put $\varphi(x)\equiv \hat a$ for all $x\in S$, for an arbitrarily fixed $\hat a\in{\bf A}$.

Let $x\in{\bf X}\setminus S$ be fixed. Since ${\bf A}=\bigcup_{i=1}^{n_1} O(a^1_i,\varepsilon_1)$ and, for each $i\in\{1,2,\ldots, n_1\}$, $\sigma^s(O(a^1_i,\varepsilon_1)|x)\in\{0,1\}$, there is the minimal index $i_1\in\{1,2,\ldots, n_1\}$ such that  \linebreak$\sigma^s(O(a^1_{i_1},\varepsilon_1)|x)=1$. We denote $\bar O^1:=O(a^1_{i_1},\varepsilon_1)$. Suppose that we have constructed a set $\bar O^k\in{\cal B}({\bf A})$ for $k=1,2,\ldots$ such that $\sigma^s(\bar O^k|x)=1$. Then we put
$$\bar O^{k+1}:=\bar O^k\cap\hat O^{k+1},$$
where $\hat O^{k+1}=O(a^{k+1}_{i_{k+1}},\varepsilon_{k+1})$ is the first one among the neighbourhoods \linebreak  $\{O(a^{k+1}_i,\varepsilon_{k+1})\}_{i=1}^{n_{k+1}}$ on which $\sigma^s(\cdot|x)$ takes the value $1$; thus $\sigma^s(\hat O^{k+1}|x)=1$. Note, $\sigma^s(\bar  O^{k+1}|x)=1$ because
\begin{eqnarray*}
1 &=& \sigma^s(\bar O^k\cup \hat O^{k+1}|x)=\sigma^s(\bar O^k|x)+\sigma^s(\hat O^{k+1}|x)-\sigma^s(\bar O^k\cap \hat O^{k+1}|x)\\
&=& 2-\sigma^s(\bar O^{k+1}|x).
\end{eqnarray*}

For the sequence $\{\bar O^k\}_{k=1}^\infty$, we have the following assertions.
\begin{itemize}
\item $\sigma^s(\bar  O^{k}|x)=1$ for all $k=1,2,\ldots$ and $\bar O^1\supseteq\bar O^2\supseteq\ldots$. Thus $\sigma^s\left(\bigcap_{k=1}^\infty \bar O^k|x\right)=1$, and hence $\bigcap_{k=1}^\infty \bar O^k\ne\emptyset$.
\item $\bigcap_{k=1}^\infty \bar O^k=\{b\}$ is a singleton because, if $b_1,b_2\in \bigcap_{k=1}^\infty \bar O^k$, then, for each $k\ge 1$, $b_1,b_2\in O(a^k_{i_k},\varepsilon_k)$ for some $i_k\in\{1,2,\ldots, n_k\}$ leading to the inequalities
$$\kappa(b_1,b_2)\le \kappa(b_1,a^k_{i_k})+\kappa(a^k_{i_k},b_2)\le 2\varepsilon_k.$$
As the result, $\kappa(b_1,b_2)\le\lim_{k\to\infty} 2\varepsilon_k=0$.
\end{itemize}

We put $\varphi(x):=b$ for that preliminarily fixed $x\in{\bf X}\setminus S$ and for $b\in{\bf A}$ such that $\bigcap_{k=1}^\infty \bar O^k=\{b\}$. As  was shown above, $\sigma^s(\{\varphi(x)\}|x)=\sigma^s(\bigcap_{k=1}^\infty \bar O^k|x)=1$; so $\sigma^s(da|x)=\delta_{\varphi(x)}(da)$ for all $x\in{\bf X}\setminus S$, that is, for $\MM^{\sigma^s}_{x_0}(dx\times{\bf A})$-almost all $x\in{\bf X}$ because $\MM^{\sigma^s}_{x_0}(S\times{\bf A})=0$. The mapping $\varphi:~{\bf X}\to{\bf A}$ is measurable because, for all $\Gamma^A\in{\cal B}({\bf A})$,
$$\{x\in{\bf A}:~\varphi(x)\in\Gamma^A\}=\left\{\begin{array}{ll}
\{x:~\sigma^s(\Gamma^A|x)=1\}\setminus S, &~ \mbox{ if } \hat a\notin\Gamma^A; \\
\{x:~\sigma^s(\Gamma^A|x)=1\}\cup S, &~ \mbox{ if } \hat a\in\Gamma^A.
\end{array}\right.$$
(Recall that the stochastic kernel $\sigma^s$ is measurable and $S\in{\cal B}({\bf X})$.)

As the result,
$$\MM(dx\times da)=\MM^{\sigma^s}_{x_0}(dx\times da)=\sigma^s(da|x)\MM^{\sigma^s}_{x_0}(dx\times{\bf A})=\delta_{\varphi(x)}(da)\MM^{\sigma^s}_{x_0}(dx\times{\bf A})$$
on  ${\cal B}({\bf X}\times{\bf A})$,
and $\MM=\MM^{\varphi}_{x_0}$ according to Lemma \ref{lB8} because the deterministic stationary strategy $\varphi$ is induced by $\MM\in{\cal D}^f$.
$\hfill\Box$

Before proving Theorem \ref{tB2}, we present several statements on  mathematical programs.

Suppose $\cal X$ is a convex compact space and $\hat{\cal C}$ is the space of $(-\infty, +\infty]$-valued bounded from below lower semicontinuous affine functions on $\cal X$. Let $R_0(\cdot),R_1(\cdot),\ldots,R_J(\cdot)\in\hat{\cal C}$ and consider the following constrained problem
\begin{eqnarray}\label{eB12p}
\mbox{Minimize over $x\in{\cal X}$: } R_0(x) \mbox{ subject to } R_j(x)\le d_j,~j=1,2,\ldots, J,
\end{eqnarray}
where $d_j\in\RR$ are fixed constants and $J\ge 1$. Here by a convex space we mean a convex subset of a cone. This definition does not involve a linear space to embed the given space in. The terms of affine functions and extreme points are understood wrt convex spaces, or say convex sets in a cone. See the definitions in \cite{archive}, where further relevant literature can be found. For all our applications here, it is sufficient to remember  that  the space of occupation measures is a convex subset of the cone of $[0,\infty]$-valued measures. When some occupation measures can take infinite values, it is difficult to embed the space of occupation measures into a convex subset of any linear space, given that we use the usual notions of addition and scalar multiplication for measures.

\begin{proposition}\label{tB1}
Consider problem (\ref{eB12p}) as was described in the above paragraph. Suppose problem (\ref{eB12p}) is non-degenerate, i.e., there is at least one point $\hat x\in{\cal X}$ satisfying all the inequalities in (\ref{eB12p}). Assume also that $\hat{\cal C}$ separates points in $\cal X$. Then there exists a solution to  problem (\ref{eB12p}) in the form $\sum_{k=1}^{J+1} \alpha_k x_k$, where $\alpha_k\in[0,1]$, $\sum_{k=1}^{J+1} \alpha_k=1$, and $x_k$ is extreme in $\cal X$ for each $k=1,2,\ldots,J+1$.
\end{proposition}
\par\noindent\textit{Proof.} See \cite[Theorem 2.1]{archive}. $\hfill\Box$

If $\mathbb E$ is a nonempty convex subset of $\RR^n$ and $u\in{\mathbb E}$, then $G(u)$ denotes the minimal face of $\mathbb E$ containing the point $u$.
A point $u\in {\mathbb E}$ is called Pareto optimal if, for each $v\in {\mathbb E}$, the componentwise inequality $v\le u$ implies that $v=u$. The collection of all Pareto optimal points is denoted by $Par({\mathbb E})$. The following result taken from \cite{FeinbergS:1996} reveals the structure of $G(u).$
\begin{proposition}\label{lb31}
Suppose ${\mathbb E}$ is a fixed nonempty convex subset of $\RR^n$, and $u\in Par({\mathbb E})$. Then the following assertions are valid.
\begin{itemize}
\item[(a)] $G(u)\subseteq Par({\mathbb E})$.
\item[(b)] For some $1\le k\le n$, there exist hyperplanes
\begin{eqnarray*}
\HH^i=\{x\in\mathbb{R}^n:~\langle x,b^i\rangle =\beta^i\},~i=1,2,\ldots, k
\end{eqnarray*}
with the following properties:
\begin{itemize}
\item[(i)] $b^i\ge 0$ for $i=1,2,\ldots, k-1$ and $b^k>0$. Here all the inequalities are componentwise.
\item[(ii)] $\HH^1$ is supporting to ${\mathbb E}^0:={\mathbb E}$ at $u$; for $i=1,2,\ldots,k-1$, ${\mathbb E}^i:={\mathbb E}^{i-1}\cap \HH^i$ and $\HH^{i+1}$ is supporting to ${\mathbb E}^i$ at $u$;
\item[(iii)] $G(u)={\mathbb E}^k:={\mathbb E}^{k-1}\cap \HH^k$.
\end{itemize}
\end{itemize}
\end{proposition}
\par\noindent\textit{Proof.}  See Lemmas 3.1 and 3.2 of \cite{FeinbergS:1996}. $\hfill\Box$

Now we are ready to present the proof of Theorem \ref{tB2}.

\par\noindent\textit{Proof of Theorem \ref{tB2}.}  The idea is to make use of Proposition \ref{tB1} and Theorem \ref{lB33}.

In the space $\cal D$ of occupation measures, we fix the topology $\rho$ as in Definition \ref{dB3}.
According to  Corollary \ref{corB3}, there is a solution $\MM^*\in{\cal D}$ to  problem (\ref{eB7}), equivalent to (\ref{ea6}), which has the form of problem (\ref{eB12p}):
\begin{itemize}
\item the space ${\cal X}={\cal D}$  is convex compact due to Proposition \ref{prB31} and Lemma \ref{lB7}(a);
 \item the mappings $R_j(\cdot)$, $j=0,1,\ldots,J$,  are non-negative, affine and lower semicontinuous  by Lemma  \ref{lB7}(b).
\end{itemize}
According to  Step 1 in the proof of \cite[Theorem 2.1]{archive}, one can accept that the point
\begin{eqnarray*}
\vec R^*:=(R_0(\MM^*),R_1(\MM^*),\ldots, R_J(\MM^*))\in\RR^{J+1}
\end{eqnarray*}
belongs to $Par({\mathbb O}\cap\RR^{J+1})$, where
\begin{eqnarray*}
{\mathbb O}:=\{\vec R(\MM)=(R_0(\MM),R_1(\MM),\ldots,R_J(\MM)),~\MM\in{\cal D}\}
\end{eqnarray*}
is the (convex) objective space.

We denote ${\mathbb E}^0={\mathbb E}:={\mathbb O}\cap\RR^{J+1}$ and emphasize that
\begin{eqnarray*}
\vec R(\MM)\in {\mathbb E}^0 \Longleftrightarrow \MM\in {\mathbb F}^0&:=&\{\MM\in{\cal D}:~\vec R(\MM)\in\RR^{J+1}\}\\
&=&\{\MM\in{\cal D}^f:~\vec R(\MM)\in\RR^{J+1}\}.
\end{eqnarray*}
The equality holds because,  due to the imposed conditions,  the component $R_{\tilde j}(\MM)$ cannot be finite if $\MM({\bf X}\times{\bf A})=+\infty$. The set ${\mathbb F}^0$, the full pre-image of ${\mathbb E}^0$ wrt the mapping $\vec R(\cdot):~{\cal D}\to (\RR\cup \{\infty\})^{J+1}$, is a face of ${\cal D}^f$: recall that the mapping $\vec R(\cdot)$ is affine.

Consider the sets ${\mathbb E}^i$, $i=0,1,\ldots, k\le J+1$  and the hyperplanes ${\mathbb H}^i=\{x\in\RR^{J+1}:~\langle x, b^i\rangle=\beta^i\}$, $i=1,2,\ldots, k$,
as in Proposition \ref{lb31} applied to $n=J+1$, ${\mathbb E}$ and $u=\vec R^*$. Let ${\mathbb F}^i$ be the full pre-image of ${\mathbb E}^i$ wrt the mapping $\vec R(\cdot)$. Note that  $\MM^*\in {\mathbb F}^i$ for all $i=0,1,\ldots,k$ because $\vec R^*=\vec R(\MM^*)\in {\mathbb E}^i$ for all $i=0,1,\ldots,k$.

Firstly, let us prove that, for each $i=0,1,\ldots,k$,   $\FF^i$,  is a (nonempty) face of ${\cal D}^f$. Roughly speaking, $\FF^{i+1}$ is a face of $\FF^i$ because ${\mathbb E}^{i+1}$ is the exposed face of ${\mathbb E}^i$.
The statement to be proved is valid for $i=0$.  Suppose it holds for some $i=0,1,\ldots, k-1$. Then
\begin{eqnarray*}
\FF^{i+1} &=& \FF^i\cap \{\MM\in{\cal D}^f:~\vec R(\MM)\in\RR^{J+1},~\langle\vec R(\MM),b^{i+1}\rangle=\beta^{i+1}\}\\
&=& \{\MM\in \FF^i:~\langle\vec R(\MM),b^{i+1}\rangle=\beta^{i+1}\}
\end{eqnarray*}
because ${\mathbb E}^{i+1}={\mathbb E}^i\cap \HH^{i+1}$.
For each $\MM\in \FF^i$, $\vec R(\MM)\in {\mathbb E}^i$, so  $\langle\vec R(\MM),b^{i+1}\rangle\ge\beta^{i+1}$ because the hyperplane $\HH^{i+1}$ supports ${\mathbb E}^i$ at $\vec R^*$. Therefore, if $\MM=\alpha \MM_1+(1-\alpha)\MM_2\in \FF^{i+1}$ for $\alpha\in(0,1)$ and $\MM_1,\MM_2\in F^i$, then
$$\langle\vec R(\MM_{1,2}) ,b^{i+1}\rangle =\beta^{i+1},\mbox{ and hence } \MM_1,\MM_2\in \FF^{i+1}.$$
Thus, $\FF^{i+1}$ is a face of  $\FF^i$ and, consequently, a face of ${\cal D}^f$ because  $\FF^i$ is a face of ${\cal D}^f$ by the induction supposition.

We have proved that $\FF^k$ is a nonempty face of ${\cal D}^f$ and $\MM^*\in \FF^k$. In fact, $\FF^k$ is the full pre-image of $G(\vec R^*)={\mathbb E}^k$, the minimal face of ${\mathbb E}={\mathbb O}\cap \RR^{J+1}$ containing $\vec R^*$: see Proposition \ref{lb31}.

Secondly, let us show that the face $\FF^k$ is closed and hence compact. Since the hyperplanes $\HH^{i+1}=\{x\in\RR^{J+1}:~\langle x, b^{i+1}\rangle=\beta^{i+1}\}$ are supporting ${\mathbb E}^i$ at $u=\vec R^*$ ($i=0,1,\ldots, k-1$), one can also write
\begin{eqnarray*}
\FF^{i+1} &=&  \{\MM\in \FF^i:~\langle\vec R(\MM),b^{i+1}\rangle\le\beta^{i+1}\}= \FF^i\cap  \{\MM\in {\cal D}:~\langle\vec R(\MM),b^{i+1}\rangle\le\beta^{i+1}\},
\end{eqnarray*}
so that
\begin{eqnarray}
\FF^k &=&  \FF^0\cap\left(\bigcap_{i=0}^{k-1} \{\MM\in {\cal D}:~\langle\vec R(\MM),b^{i+1}\rangle\le\beta^{i+1}\}\right)\nonumber\\
&=& \bar \FF^0 \cap\left(\bigcap_{i=0}^{k-2} \{\MM\in {\cal D}:~\langle\vec R(\MM),b^{i+1}\rangle\le\beta^{i+1}\}\right),\label{eB99}
\end{eqnarray}
where
\begin{eqnarray*}
\bar \FF^0 &:=& \FF^0\cap \{\MM\in{\cal D}:~\langle\vec R(\MM),b^{k}\rangle\le\beta^{k}\}= \{\MM\in{\cal D}:~\langle\vec R(\MM),b^{k}\rangle\le\beta^{k}\}.
\end{eqnarray*}
The second equality  holds because $\vec R(\MM)\ge 0$ and $b^k>0$: if $\langle\vec R(\MM),b^k\rangle\le \beta^k$, then $\vec R(\MM)\in\RR^{J+1}$;   so that  $\MM\in \FF^0$.
Note, $\bar \FF^0$ is not necessarily a face of $\cal D$. The space $({\cal D},\rho)$ is compact and all the mappings $\langle\vec R(\cdot),b^{i+1}\rangle:~{\cal D}\to\bar\RR^0_+$, $i=0,1,\ldots, k-1$, are lower semicontinuous by Lemma \ref{lB7}. Therefore, the set $\bar \FF^0$ is closed, and  the face $\FF^k$ of ${\cal D}^f$ is closed by (\ref{eB99}), hence, compact, as the closed subset of the compact $\cal D$. %

The space $\hat{\cal C}$ of $(-\infty,+\infty]$-valued bounded from below lower semicontinuous affine functions on ${\cal D}^f$ separates points in ${\cal D}^f$. Indeed,
if $\MM_1\ne\MM_2$ are two measures from ${\cal D}^f$, then
$$C(\MM_1):=\int_{{\bf X}\times{\bf A}} c(x,a)\MM_1(dx\times da)\ne \int_{{\bf X}\times{\bf A}} c(x,a)\MM_2(dx\times da)=:C(\MM_2)$$
for some non-negative bounded continuous function $c(\cdot,\cdot):~{\bf X}\times{\bf A}\to\RR$. (See Lemma 2.3 of \cite{Varadarajan:1958}, Theorem 5.9 of \cite{par} and Proposition 7.18 of \cite{Bertsekas:1978}.)
The mapping $C(\cdot)$ is non-negative, lower semicontinuous and affine by Lemma \ref{lB7}(b), so that the desired assertion  follows.

Since the compact face $\FF^k\subseteq{\cal D}^f$ contains $\MM^*$, one can consider problem  (\ref{eB7}), on $\FF^k$, not on ${\cal D}$. This problem
\begin{eqnarray}
R_0(\MM) &:=& \int_{{\bf X}\times{\bf A}}r_0(x,a)\MM(dx\times da)\to \min_{\MM\in \FF^k} \label{eB72}\\
\mbox{s.t.}~R_j(\MM) &:=& \int_{{\bf X}\times{\bf A}}r_j(x,a)\MM(dx\times da)\le d_j,~~j=1,2,\ldots, J,\nonumber
\end{eqnarray}
satisfies all the conditions of Proposition \ref{tB1}:
\begin{itemize}
\item the space $\FF^k$ is convex compact;
\item the mappings $R_j(\cdot):~\FF^k\to\bar\RR^0_+$, $j=0,1,\ldots, J$, are non-negative, lower semicontinuous by Lemma \ref{lB7}(b), and affine;
\item problem (\ref{eB72}) is non-degenerate, as the measure $\MM^*\in \FF^k$ satisfies all the constraints;
\item The space $\hat{\cal C}$ separates points in $\FF^k\subseteq{\cal D}^f$.
\end{itemize}

According to Proposition \ref{tB1}, there exists a solution to  problem (\ref{eB72}) (hence, to  problem    (\ref{eB7})) in the form $\sum_{l=1}^{J+1} \alpha_l\MM_l$, where $\alpha_l\in[0,1]$, $\sum_{l=1}^{J+1}\alpha_l=1$, and $\MM_l$ is extreme in $\FF^k$ for each $l=1,2,\ldots,J+1$. Since $\FF^k$ is a face of ${\cal D}^f$, $\MM_l$ is extreme also in ${\cal D}^f$ for each $l=1,2,\ldots,J+1$ and equals $\MM^{\varphi_l}_{x_0}$ for some deterministic stationary strategy $\varphi_l$ in accordance with Theorem \ref{lB33}. Therefore, the mixture $\hat\MM:=\sum_{l=1}^{J+1}\alpha_l\MM^{\varphi_l}_{x_0}$ solves  problem  (\ref{eB7}), and the corresponding strategic measure $\hat\PP:=\sum_{l=1}^{J+1}\alpha_l\PP^{\varphi_l}_{x_0}$ solves  problem (\ref{ea6}). $\hfill\Box$

\end{document}